%% file: intro.tex
\documentclass[reqno]{amsart}
\usepackage{pictex}
\usepackage{mathrsfs}
\usepackage{hyperref}
\usepackage{color}
\usepackage{graphicx}

\begin{document}

% #############################################
%
%           Macros, etc
%
% #############################################

\def\labelenumi{(\theenumi)}

\newtheorem{thm}{Theorem}[section]
\newtheorem{lem}[thm]{Lemma}
\newtheorem{conj}[thm]{Conjecture}
\newtheorem{cor}[thm]{Corollary}
\newtheorem{add}[thm]{Addendum}
\newtheorem{prop}[thm]{Proposition}
\theoremstyle{definition}
\newtheorem{defn}[thm]{Definition}
\theoremstyle{remark}
\newtheorem{rmk}[thm]{Remark}
\newtheorem{example}[thm]{{\bf Example}}
\newtheorem{Question}[thm]{{ \bf Question}}

\newcommand{\SurfG}{\Sigma_g}
\newcommand{\SurfGN}{\Sigma_{g,n}}
\newcommand{\TriangG}{T_g}
\newcommand{\TriangGOne}{T_{g,1}}
\newcommand{\ProjG}{\mathcal{P}_g}
\newcommand{\TeichG}{\mathcal{T}_g}
\newcommand{\CirclePackGTau}{\mathsf{CPS}_{g,\tau}}
\newcommand{\CrossRatio}{{\bf c}}
\newcommand{\CrossRatioGTau}{\mathcal{C}_{g,\tau}}
\newcommand{\CrossRatioOneTau}{\mathcal{C}_{1,\tau}}
\newcommand{\DeformGTau}{\mathcal{C}_{g,\tau}}
\newcommand{\Forget}{\mathit{forg}}
\newcommand{\Uniform}{\mathit{u}}
\newcommand{\Section}{\mathit{sect}}
\newcommand{\SLTwoC}{\mathrm{SL}(2,{\mathbb C})}
\newcommand{\SLTwoR}{\mathrm{SL}(2,{\mathbb R})}
\newcommand{\SUTwo}{\mathrm{SU}(2)}
\newcommand{\PSLTwoC}{\mathrm{PSL}(2,{\mathbb C})}
\newcommand{\GLTwoZ}{\mathrm{GL}(2,{\mathbb Z})}
\newcommand{\PGLTwoZ}{\mathrm{PGL}(2,{\mathbb Z})}
\newcommand{\GLTwoC}{\mathrm{GL}(2,{\mathbb C})}
\newcommand{\PSLTwoR}{\mathrm{PSL}(2,{\mathbb R})}
\newcommand{\PGLTwoR}{\mathrm{PGL}(2,{\mathbb R})}
\newcommand{\GLTwoR}{\mathrm{GL}(2,{\mathbb R})}
\newcommand{\PSLTwoZ}{\mathrm{PSL}(2,{\mathbb Z})}
\newcommand{\SLTwoZ}{\mathrm{SL}(2,{\mathbb Z})}
\newcommand{\nnn}{\noindent}
\newcommand{\MCG}{{\mathcal {MCG}}}
\newcommand{\MMap}{{\bf \Phi}_{\mu}}
\newcommand{\HH}{{\mathbb H}^2}
\newcommand{\TT}{{T_{1,0}}}
\newcommand{\X}{{\mathcal  X}}
\newcommand{\B}{{\mathcal  B}}
\newcommand{\E}{{\mathcal  E}}
\newcommand{\C}{{\mathscr C}}
\newcommand{\T}{{\mathscr T}}
\newcommand{\M}{{\mathscr M}}
\newcommand{\F}{{\mathscr F}}
\newcommand{\CC}{{\mathbb C}}
\newcommand{\RR}{{\mathbb R}}
\newcommand{\Li}{{{\mathrm{Li}_2}}}

\renewcommand{\L}{{\mathcal  L}}
\newcommand{\G}{{\mathcal  G}}
\newcommand{\R}{{\mathcal  R}}
\newcommand{\Q}{{\mathbb Q}}
\newcommand{\ZZ}{{\mathbb Z}}
\newcommand{\PL}{{\mathscr {PL}}}
\newcommand{\GP}{{\mathcal {GP}}}
\newcommand{\GT}{{\mathcal {GT}}}
\newcommand{\GQ}{{\mathcal {GQ}}}
\newcommand{\EE}{{{\mathcal E}(\rho)}}
\newcommand{\HHH}{{\mathbb H}^3}
\newcommand{\tr}{{\rm tr\, }}

\def\square{\hfill${\vcenter{\vbox{\hrule height.4pt \hbox{\vrule width.4pt
height7pt \kern7pt \vrule width.4pt} \hrule height.4pt}}}$}

\newenvironment{pf}{\noindent {\sl Proof.}\quad}{\square \vskip 12pt}

\title{A dilogarithm identity on moduli spaces of curves}
%for (closed) hyperbolic surfaces]%
%{A new identity for closed hyperbolic surfaces}
\author[F. Luo \& S. P. Tan]{Feng Luo \& Ser Peow Tan
 }
\address{Department of Mathematics,\\
Rutgers  University,\\
Piscataway, NJ 08854, USA}
\email{fluo@math.rutgers.edu}%
\address{Department of Mathematics \\ National University of Singapore \\
2 Science Drive 2 \\ Singapore 117543}
\email{mattansp@nus.edu.sg}%

\dedicatory{To Michael Freedman on the occasion of his sixtieth
birthday}



\thanks{The first author  is
partially supported by the NSF (NSF-CCF-0830572) of USA.  The
second author is partially supported by the National University of
Singapore academic research grant R-146-000-133-112}

% #############################################
%
%                  Abstract
%
% #############################################

\begin{abstract}
We establish an identity for compact hyperbolic surfaces with or
without boundary whose terms depend on the dilogarithms of the
lengths of simple closed geodesics in all 3-holed spheres and
1-holed tori in the surface.

\end{abstract}

\maketitle

\input sec1.tex

\input sec2.tex

\input sec3.tex

\input sec4.tex

%\input sec5.tex

\input sec6.tex

\input ref.tex

\input appendix.tex

\end{document}

%% file: sec1.tex
\section{introduction}

\subsection{Statement of  results}
In \cite{mcshane1}, McShane established a remarkable identity for
the lengths of simple closed geodesics in hyperbolic surfaces with
cusp ends. Since then there have been many generalizations of
McShane's identity, for example, to hyperbolic surfaces with
geodesic boundaries \cite{mirz}, \cite{tan1} and cone
singularities \cite{tan1}. Mirzakhani also found fantastic
applications of these identities to the computation of the volumes
of moduli spaces of bordered Riemann surfaces.  There has been
much research since then towards finding a McShane type identity
for closed hyperbolic surfaces. In \cite{mcshane2} and
\cite{tan1}, McShane and Tan et al established such an identity for
closed hyperbolic surfaces of genus 2. However, the techniques
used there do not generalize as they depend crucially on the fact
that every genus 2 surface admits a hyperelliptic involution. The
goal of this paper is to establish a McShane type identity for
simple closed geodesics on \emph{any} closed hyperbolic surface.
Our result for the genus 2 case is different from that given in
\cite{mcshane2} or \cite{tan1}. The generalization of our identity
to surfaces with cusps or geodesic boundary also differ from those
in \cite{mcshane1} or \cite{mirz}. This seems to suggest that
there are possibilities of producing many different McShane type
identities for  hyperbolic surfaces. We expect that the identity
found here will have applications towards the study of the moduli
space of curves.

The identity that we produce involves the dilogarithm of the
lengths of simple closed geodesics in all 1-holed tori and 3-holed
spheres in the surface. Our work is motivated by \cite{mcshane1},
\cite{mirz}, \cite{tan1} and \cite{bridge}.  In \cite{bridge},
Bridgeman considers compact hyperbolic surfaces with non-empty
geodesic boundary and geodesic paths starting and ending at the
boundary. Our approach is similar to that of \cite{bridge} in two
aspects. First, we consider the unit tangent bundle instead of the
surface itself, and second the identity obtained involves
dilogarithm functions. In fact, we use Bridgeman's work in
producing the main identity. The main idea in arriving at our
identity is also closely related to the interpretation and proof
of McShane's identity in \cite{tan1}.

In this paper, we consider oriented surfaces.  For a hyperbolic
surface $F$, a compact embedded subsurface $\Sigma \subset F$ is
said to be \emph{ geometric }if the boundaries of $\Sigma$ are
geodesic and \emph{proper} if the inclusion map $i: \Sigma
\rightarrow F$ is injective. Furthermore call a surface \it simple
\rm if it is a 3-holed sphere or one hole torus (both of them have
Euler characteristic $-1$).  Our main result is the following:

\begin{thm}\label{thm:main}
Let $F$ be a closed hyperbolic surface of genus $g \ge 2$. There
exist functions $f$ and $g$ involving the dilogarithm of the
lengths of the  simple closed geodesics  in a 3-holed sphere or 1-holed torus, such that

\begin{equation}\label{eqn:maintheorem}
   \sum_{P} f(P)+\sum_{T} g(T)= 8\pi^2(g-1)
\end{equation}
where the first sum is over all properly embedded geometric 3-holed
spheres $P \subset F$, the second sum is over all properly embedded
geometric 1-holed tori $T\subset F$.

\end{thm}
The definitions of the functions $f$ and $g$ in the identities are
given in \S \ref{s:functions}. The right-hand-side in (1) is the
volume of the unit tangent bundle over the surface $F$.

\bigskip

\noindent Remarks: \begin{enumerate}

                     \item Each $T$ in the second summand can be cut along
                      simple closed geodesics into a 3-holed sphere. These 3-holed spheres
                      are not properly embedded.
                     \item Bridgman's identity \cite{bridge} does not extend to
                     closed hyperbolic surfaces without boundary. Nonetheless,
                     the terms involved in our identity are similar to those of Bridgeman's
                     in the sense that they involve the Roger's dilogarithm function.
                     \item Our identity can be thought of as a hybrid of both the McShane
                     and Bridgeman identities.
                    It can also be thought of as an identity on the moduli space $\M_g$ of curves,
                    rather than the Teichm\"uller space $\T_g$,
                    as the mapping class group has the effect of permuting the terms
                    in the summands.
                    \item  The theorem can be extended to
                    hyperbolic surfaces with geodesic boundary and
                    cusp ends. The expression is more complicated
                    though. See theorem 1.2 below.
                    \item We have been informed by G. McShane that he
                    and D. Calegari have recently obtained results similar to theorem 1.1.
                   \end{enumerate}

\bigskip

\subsection{Basic idea of the proof}
The key idea is to decompose the unit tangent bundle $S(F)$ of a
closed hyperbolic surface $F$ according to, and indexed by, the
properly embedded geometric 1-holed tori and 3-holed spheres in
$F$. The decomposition is measure theoretic in the sense that we
will ignore a measure zero set in $S(F)$. Here is the way to
produce the decomposition. For a unit tangent vector $v \in S(F)$,
consider the unit speed geodesic rays $g_v^+(t)$ and $g_v^-(t)$ $(
t \geq 0)$ determined by $\pm v$. If the vector $v$ is generic,
then both rays will self intersect transversely by the ergodicity
of the geodesic flow. This vector $v$ will determine a canonical
graph $G(v)$ as follows. Consider the path $A_t=g_v^-([0,t]) \cup
g_v^+([0,t])$ for $t
>0$ obtained by letting the geodesic rays $g_v^{-}$ and
$g_v^+$ grow at equal speed from time $0$ to $t$. Let $t_1>0$ be
the smallest positive number so that $A_{t_1}$ is not a simple
arc. Say $g_v^+(t_1) \in g_v^-([0,t_1]) \cup g_v^+([0, t_1))$.
Next, let  $t_2 \geq t_1$ be the next smallest time so that
$g_v^-(t_2) \in g_v^-([0, t_2)) \cup g_v^+[0, t_1])$.

\begin{figure}[ht!]
\centering
\includegraphics[scale=0.55]{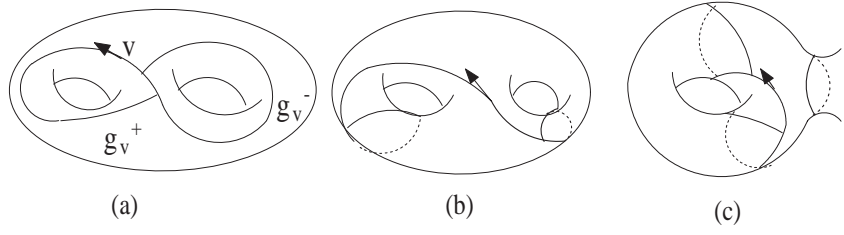}
\caption{creation of spine} \label{figure 1}
\end{figure}

%%The description given below is not quite accurate so I replaced with that above.
%Let $t_1$ be the smallest time so that $\gamma_+$ or $\gamma_-$, say
%$\gamma_-$ becomes self-intersect. Next, extending $\gamma_+(t)$
%from $[0, t_1]$ further to the first time $t_2 \geq t_1$ so that
%$g(t_2)$ is in $\gamma_-([0,t_1]) \cup \gamma_+([0,t_2))$.
The union $g_v^-[0, t_2] \cup g_v^+[0, t_1]$ is the graph, denoted
by $G(v)$ associated to $v$. Its Euler characteristic is $-1$. The
graph $G(v)$ is contained in a unique properly embedded geometric
subsurface $\Sigma(v)$ which is either a 1-holed torus or a
3-holed sphere in $F$. Furthermore either the graph $G(v)$ is a
deformation retract of $\Sigma(v)$, or $\Sigma(v)$ is a 1-holed
torus so that $\Sigma(v) -G(v)$ is a union of two annuli (figure
1(c)). By abuse of notation, we will say in this case that $G(v)$
is also a spine for $\Sigma(v)$. In this way, we produce a
decomposition of the unit tangent bundle $S(F)$. Namely,
generically, each vector $v \in S(F)$ is in a unique geometric
1-holed torus $T$ or a 3-holed sphere $P$ so that $G(v)$ is a
spine for the subsurface. It remains to calculate for a simple
hyperbolic surface $\Sigma$ the volume of the set of all unit
tangent vectors $v$ in $\Sigma$ so that $G(v)$ is a spine for
$\Sigma$. It turns out the volume of this set can be explicitly
calculated using the dilogarithm and the lengths of simple closed
geodesics in $\Sigma$.

\bigskip
\subsection{Extension to non-closed hyperbolic surfaces}

The identity can be extended to finite area hyperbolic surfaces
$F_{g,n}^r$ of genus $g$ with $n$ geodesic boundary components and
$r$ cusps, by modifying the function $f$ for 3-holed spheres $P$
whose boundaries become peripheral.
%The main idea is that if $G(v)$
%intersects the boundary of $P$ in components which are contained
%in the boundary of $F_{g,n}^0$, we should calculate the volume of
%all the unit vectors $v$ such that the union of $G(v)$ with these
%boundary components form a spine for $P$.

Let $\F_{g,n}^r$ denote the set of all marked hyperbolic
structures a surface $\Sigma_{g,n}^r$ of genus $g$ with $n$
boundary components and $r$ punctures so that the boundaries are
geodesics and punctures are cusps.  Let $F_{g,n} = F^0_{g,n}$ and
$\F_{g,n} = \F_{g,n}^0$.

We have:

\begin{thm}\label{thm:mainwithboundary}
Let $F \in \F_{g, n}^0$ be a hyperbolic surface with geodesic
boundaries so that its Euler characteristic is strictly less than
$-1$. There exist functions ${\hat f}, {\bar f}, f:\F_{0,3}
\rightarrow \RR_+$, $g:\F_{1,1} \rightarrow \RR_+$ such that

\begin{equation}\label{eqn:maintheorem2}
  \sum_{\hat P} {\hat f}(\hat P)+
  \sum_{\bar P} {\bar f}(\bar P)+\sum_{P} f(P)+\sum_{T} g(T)=
  4\pi^2(2g-2+n)
\end{equation}
where the first sum is over all properly embedded geometric pairs
of pants $\hat P \subset F$ with exactly one boundary component in
 $\partial F$, the second is over all properly embedded geometric
pairs of pants $\bar P \subset F$ with exactly two boundary
components in $\partial F$, the third sum is over all properly
embedded geometric pairs of pants $P \subset F$ such that
$\partial  P \cap
\partial F=\emptyset$, the fourth sum is over all properly embedded
geometric one holed tori $T \subset F$.

Furthermore, if lengths of $k$ boundary components of $F_{g,n}$
tend to zero, then each term and each summation in (2) converge.
The limit is the identity for all hyperbolic surfaces
$F_{g,n-k}^k$ of genus $g$ with $n-k$ geodesic boundary and $k$
cusps.
\end{thm}

The right hand side of (2) is the volume of the unit tangent
bundle over $F$.

\bigskip
\subsection{Plan of the paper} In section \ref{s:functions}, we
define the functions $f$, $g$, $\hat f$, $\bar f$   in
(\ref{eqn:maintheorem}) and (\ref{eqn:maintheorem2}). In section
\ref{s:decomposing}, we describe how to decompose the unit tangent
bundle $S(F)$ of the surface $F$ by showing how each $v \in S(F)$
generates a spine for a simple subsurface $\Sigma \subset F$.  In
section \ref{s:measureofdecompostion}, for simple subsurfaces
$\Sigma \subset F$, we identify the subset of the unit tangent
vectors in $S(\Sigma)$ which generate spines for $\Sigma$ with a
subset of $S(\HH)$. In section 5, we derive the formula for the
measure of the set studied in section 4, thereby giving the
formulas
 for $f$, $g$,  $\hat f$ and $\bar f$. Finally, in the appendix,
we first give an interpretation of the pentagon relation for the
dilogarithm function in terms of lengths of right-angled
hyperbolic pentagons and then explain why different rules for
generating the spines $G(v)$ will result in the same theorem
\ref{thm:main}.

% we
%talk about two sets of independent results related to the main
%theorem. The first is a new dilogarithm pentagon relation
%involving the side lengths of a right angled hyperbolic pentagon,
%the second is an exposition of how different rules for generating
%the spines $G(v)$ may produce different decompositions of the unit
%tangent bundle $S(F)$, but which nonetheless give the same
%identities in Theorems \ref{thm:main} and
%\ref{thm:mainwithboundary}.
%, and finally we give a brief description of the work in \cite{HuTan} which gives an identity for the 1-holed torus %which allows us to rewrite $g(T)$ as an infinite sum over all non-boundary parallel simple closed geodesics on $T$.

\bigskip
\noindent {\it Acknowledgements.} The authors thank Norman Do,
Greg McShane, Yasushi Yamashita and Ying Zhang for helpful
conversations. They also thank Bill Goldman, Scott Wolpert and
Danny Calegari for helpful comments on the first version of the
paper. Part of this work was carried out during the program on the
Geometry, Topology and Dynamics of Character Varieties held in the
Institute for Mathematical Sciences in NUS in the summer of 2010,
the authors thank the IMS for its support. The first author
thanks the Center of Mathematical Science at Qinghua University,
China where part of the work was carried out.

%\end{document}

%% file: sec2.tex
\section{Definitions of the functions}\label{s:functions}
In this section we define $f$, $g$, $\hat f$, $\bar f$   in the
identities (\ref{eqn:maintheorem}) and (\ref{eqn:maintheorem2}).

\subsection{Dilogarithm and Roger's dilogarithm functions}\label{ss:dilog}
We first recall the dilogarithm function $\Li$ and the Roger's
dilogarithm function $\L$. See \cite{lewin} for more details.

The dilogarithm function $\Li$ is defined for $|z|<1$ by the Taylor series
\begin{equation}\label{eqn:dilog1}
    \Li(z)=\sum_{n=1}^{\infty} \frac{z^n}{n^2}
\end{equation}
so that for $x \in \bold R$ with $x<1$,
%$\Li(z)+\Li(-z)=\frac{1}{2}\Li(z^2). $
$   \Li(x)=-\int_{0}^x\frac{\log (1-z)}{z}dz.$

The Rogers $\L$-function is defined by
\begin{equation}\label{eqn:Rogersdilog}
    \L (x)=\Li (x)+\frac{1}{2} \log(|x|)\log (1-x)
\end{equation}
so that $    \L'(z)=-\frac{1}{2}\left( \frac{\log
(1-z)}{z}+\frac{\log (z)}{1-z} \right)$ and $\L(0)=0$. It
satisfies $\L(x) + \L(1-x) =\pi^2/6$ for $0<x<1$. The fundamental
identity which characterizes the function $\L(x)$ is the following
pentagon relation, for $x,y \in (0,1)$,

\begin{equation}\label{eqn:pentagon}
\L(x)+\L(y) +\L(1-xy) + \L(\frac{1-x}{1-xy}) +\L(\frac{1-y}{1-xy})
 =\frac{\pi^2}{2}.
\end{equation}
A geometric interpretation of (\ref{eqn:pentagon})  in terms of the lengths of
right-angle hexagon is given in the appendix.

%Some basic identities for $\L(x)$ are given by:

%\begin{eqnarray}\label{eqn:Rogersrelation1}
%% \nonumber to remove numbering (before each equation)
%  \L(0) &=& 0 \\
%  \L(1) &=& \frac{\pi^2}{6} \\
% \L(x)+\L(1-x) &=& \frac{\pi^2}{6} \qquad 0 \le x \le 1 \\
%  \L(-x)+\L(-x^{-1}) &=& -\frac{\pi^2}{6} \qquad  x >0 \\
% \L\left(\frac{-x}{1-x} \right) &=& -\L(x) \qquad 0<x<1 \\
% \Longrightarrow \L\left (\frac{1}{\cosh ^2 (x)}\right) &=& -\L \left(\frac{-1}{\sinh^2(x)}\right)\\
% \L(x)+\L(y) &=& \L(xy)+ \L\left(\frac{x(1-y)}{1-xy} \right)+\L\left(\frac{y(1-x)}{1-xy} \right)\\
%\frac{\pi^2}{3} &=&\L(xy)+\L(1-y)+\L(1-x)+\L\left(\frac{x(1-y)}{1-xy} \right) + \L\left(\frac{y(1-x)}{1-xy} \right) \\
%\end{eqnarray}

%We also have the pentagonal relations:
%\begin{eqnarray}\label{eqn:Rogersrelation2}
%% \nonumber to remove numbering (before each equation)
% \L(xy)+\L(1-y)+\L(1-x)+\L\left(\frac{x(1-y)}{1-xy} \right) + \L\left(\frac{y(1-x)}{1-xy} \right) &=& \frac{\pi^2}{3}  \\
% \L(x)+\L(1-y)+\L\left(\frac{y-x}{y} \right) + \L\left(\frac{y-x}{1-x} \right)+ \L\left(\frac{x(1-y)}{y(1-x)} \right) &=& \frac{\pi^2}{3}
%\end{eqnarray}

\subsection{Length invariants of 3-holed spheres}\label{ss:pop}
\begin{figure}[ht!]
\centering
\includegraphics[scale=0.5]{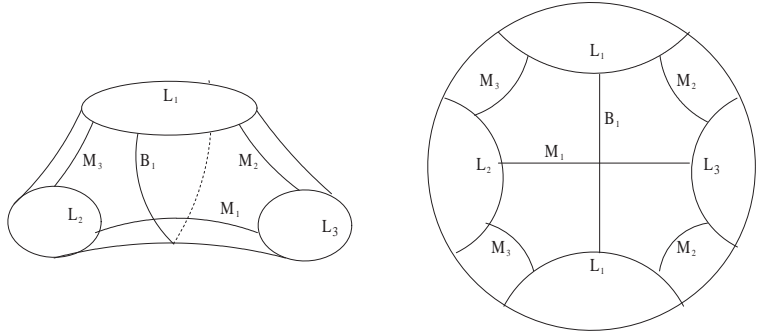}
\caption{3-holed sphere} \label{figure 2}
\end{figure}

Let $P\in \F_{0,3}$ be a hyperbolic 3-holed sphere with geodesic
boundaries $L_1, L_2, L_3$. For $\{i,j,k \}=\{1,2,3\}$,  let $M_i$
be the shortest geodesic arc between $L_j$ and $L_k$, and $B_i$
the shortest non-trivial geodesic arc from $L_i$ to itself. Note
that $M_i$ and $B_i$ are orthogonal to $\partial P$. See figure \ref{figure 2}.
We define:
\begin{itemize}
  \item $l_i$ to be the length of  $L_i$.
  \item $m_i$ to be the length of $M_i$.
  \item $p_i$ to be the length of $B_i$.
\end{itemize}

%{\bf Feng's remark: I am  replacing $P_i$ by $B_i$ since $P_i$ may
%be confusing with the pants $P$ itself.  I will retain $p_i$.}
%{\bf Tan's remark:Ok, I've also changed to 3-holed sphere to follow the intro. However, the labelling of the figures have not been changed, namely $B_i$ is still labeled $P_i$}

 Note
that $P$ is decomposed into two right-angled hyperbolic hexagons
with cyclically ordered side-lengths
 \{$\frac{l_1}{2}, m_3, \frac{l_2}{2}, m_1, \frac{l_3}{2}, m_2\}$ by cutting along the $M_i$.
 Furthermore, cutting along each $B_i$ decomposes the two
  hexagons into 2 right-angled pentagons. See figure \ref{figure 2}.
%   in three J
% different ways, see figure 2.

The sine and cosine rules for right angled hexagons and pentagons
say that for $\{i,j,k \}=\{1,2,3\}$,

\begin{equation}\label{eqn:sinerule}
    \frac{\sinh m_i}{\sinh (l_i/2)}=\frac{\sinh m_j}{\sinh (l_j/2)}=\frac{\sinh m_k}{\sinh (l_k/2)}
\end{equation}
\begin{eqnarray}\label{eqn:cosinerule}
  %\cosh(l_i/2)\sinh m_j\sinh m_k &=& \cosh m_i+\cosh m_j\cosh m_k %\\
  \cosh m_i \sinh (l_j/2) \sinh (l_k/2) &=& \cosh (l_i/2)+\cosh (l_j/2)\cosh (l_k/2)
\end{eqnarray}

\begin{equation}\label{eqn:perpendicular}
   \cosh (p_k/2)=\sinh(l_i/2)\sinh m_j
\end{equation}

In particular, all lengths $m_i, p_i$ can be expressed in terms of
$l_1,l_2$ and $l_3$.

%{\bf Tan's remark: I added the sentence above to
%emphasize that all functions really are functions of the lengths of the 3 geodesic loops} OK

\subsection{Length invariants of 1-holed tori}\label{ss:OHT}
\begin{figure}[ht!]
\centering
\includegraphics[scale=0.5]{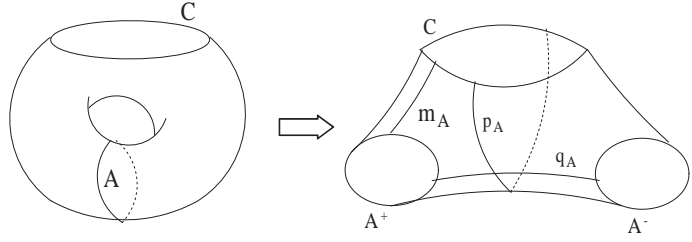}
\caption{cutting 1-holed torus into a 3-holed sphere}
\label{figure 3}
\end{figure}

Let  $T\in \F_{1,1}$ be a hyperbolic 1-holed torus with  boundary
component $C$. For any non-boundary parallel simple closed
geodesic $A$ on $T$, cutting $T$ along $A$ gives a hyperbolic pair
of pants $P_A$ with boundary geodesics $C, A^+$ and $A^-$, see
figure \ref{figure 3}. Let
\begin{itemize}
  \item $c$  be the length of $C$
  \item $a$  be the length of $A$
  \item $m_{A}$ be the shortest distance between $C$ and $A^+$ in $P_A$ (or $A^-$)
  \item $p_{A}$  be the length of the shortest non-trivial geodesic arc from $C$ to  $C$ in $P_A$
    \item $q_A$  be the length of the shortest non-trivial path from $A^+$ to $A^-$ in $P_A$.
\end{itemize}

%{\bf Tan's remark: $m_A$ and $p_A$ are lower case, they are labeled in the figure as upper case}

\subsection{The main functions and the identity}\label{ss:identity}

The functions $f$ and $g$ in theorem 1.1 are given as follows.
%{\bf Feng: I have removed the function b(x) to avoid introducing
%too many notations. It seems $\L(\frac{1}{\cosh^2(x/2)})$ is
%simple enough compare with f,g's}
%{\bf Tan's remark: Yes, I sort of agree, what do you think of
%adopting the shorthand notation, with some explanation beforehand
%of $s(x)$, $c(x)$, $t(x)$ where $s,c$ and $t$ are the hyperbolic
%sine, cosine and tangent functions respectively. Seems it would
%make the formulas neater, although may introduce some confusion?
%F: I will leave them as they are since we are not use these
%formulas extensively in the paper.}
%best expressed in
%terms of two functions which are simpler:
% We first define (following Bridgeman \cite{bridgeman}):
%\begin{equation}\label{eqn:bridgean}
 %   b(x):=\L\left(\frac{1}{\cosh^2 \frac{x}{2}}\right)
%\end{equation}
We first define the lasso function $La(l,m)$ to be
\begin{equation}\label{eqn:lasso}
La(l,m)=\L(y) - \L(\frac{1-x}{1-xy}) + \L(\frac{1-y}{1-xy})
\end{equation}
where $x =e^{-l}$ and $y=\tanh^2(m/2)$. The three terms  above
appear in (\ref{eqn:pentagon}).
%\begin{eqnarray*}\label{eqn:lasso}
%\L(1- e^{-l} \tanh^2(\frac{m}{2})) +2\L( (e^{-l}-1)
%\sinh^2(\frac{m}{2}) ) -\L(1-e^{-l}) -2\L(-\sinh^2(\frac{m}{2}) ).
%\end{eqnarray*}
%\begin{eqnarray*}\label{eqn:lasso}
%    &&La(l,m) \\
%    &:=& \L(1-e^{-l}\tanh^2(\frac{m}{2}))+2\L( (e^{-l}-1)\sinh^2(\frac{m}{2}))+\L(e^{-l}) -2\L (-\sinh^2(\frac{m}{2}))-\frac{\pi^2}{6}\\
%    &=& \L(1- e^{-l} \tanh^2(\frac{m}{2})) +2\L( (e^{-1}-1) \sinh^2(\frac{m}{2})
 %   ) -\L(1-e^{-l})
%-2\L(-\sinh^2(\frac{m}{2}) )\\
%    &=& \L(1-e^{-l}\tanh^2(\frac{m}{2}))+2\L( (e^{-l}-1)\sinh^2(\frac{m}{2}))+
 %   \L(e^{-l})-\L (-\sinh^2(\frac{m}{2}))-
 %   \L(\frac{1}{\cosh^2(\frac{m}{2})})\\
 %   &=& \L(1-e^{-l}\tanh^2(\frac{m}{2}))+2\L(
 %   (e^{-l}-1)\sinh^2(\frac{m}{2}))+\L(e^{-l})+\frac{\pi^2}{6}-2\L(\frac{1}{\cosh^2(\frac{m}{2})})\\
 %   &=& \L(e^{-l})-\L(e^{-l}\tanh^2(\frac{m}{2}))+2\L( (e^{-l}-1)\sinh^2(\frac{m}{2})) -2\L (-\sinh^2(\frac{m}{2}))
%\end{eqnarray*}

Now for $P \in \F_{0,3}$ with length invariants $l_i$, $m_i$ and
$p_i$, as given in \S \ref{ss:pop},  and $1 \le i,j \le 3$, we
define
\begin{equation}\label{eqn:definitionofffirst}
    f(P):=4\pi^2 - 8 \left\{ \sum_{i=1}^3\big( \L(\frac{1}{\cosh^2(m_i/2)})+ \L(\frac{1}{\cosh^2(p_i/2)})\big)
    +\sum_{i\neq j}La(l_i,m_j)\right\}
    \end{equation}
    \begin{equation}\label{eqn:definitionoff}
 =8 \left[\sum_{ i \neq j} (\L(\frac{1-x_i}{1-x_i y_j})
-\L(\frac{1-y_j}{1-x_i y_j}))
-\sum_{k=1}^3(\L(y_k)+\L(\frac{1}{\cosh^2(p_k/2)}))\right]
\end{equation}
\begin{equation}\label{eqn:deffp}
=4 \sum_{i \neq j}[2\L(\frac{1-x_i}{1-x_i y_j})
-2\L(\frac{1-y_j}{1-x_i y_j})-\L(y_j)
-\L(\frac{(1-y_j)^2 x_i}{(1-x_i)^2y_j})]
\end{equation}
where  $x_i = e^{-l_i}$, $y_i=\tanh^2(m_i/2)$ and by (\ref{eqn:perpendicular}), $
\frac{1}{\cosh^2(p_k/2)}=\frac{(1-y_j)^2 x_i}{(1-x_i)^2y_j}$ for
$\{i,j,k\}=\{1,2,3\}$.

 For $T \in \F_{1,1}$ with boundary geodesic $C$,
we define
\begin{equation}\label{eqn:definitionofg}
    g(T):=4\pi^2-8\sum_{A} ( \L(\frac{1}{\cosh^2(p_{A}/2)})+ 2 La(a,m_{A}) \big)
\end{equation}
where the sum is taken over all non-boundary parallel simple
closed geodesics $A$ on $T$ and $c$, $a$, $p_{A}$ and $m_{A}$ are defined as
in
\S\ref{ss:OHT}. 
A further simplification of $g(T)$ is obtained recently, see
\cite{HuTan} for details.
\begin{thm} (\cite{HuTan}) \label{thm:simplificationofGT}
\begin{eqnarray}\label{eqn:simplificationofgT}
\nonumber  g(T)&=&\sum_A\big\{4\pi^2-8\big[2\L(\frac{1}{\cosh^2(m_A/2)})+\L(\frac{1}{\cosh^2(q_A/2)})+\L(\frac{1}{\cosh^2(p_{A}/2)})\\
    && \qquad \qquad+2La(c/2,m_A)+ 2 La(a,m_{A})\big]\big\}
\end{eqnarray}
where the sum is taken over all non-boundary parallel simple
closed geodesics $A$ on $T$, and $c$, $a$, $p_{A}$, $m_{A}$
 and $q_A$ are defined in \S\ref{ss:OHT}.
\end{thm}

Identities (\ref{eqn:definitionoff}),(\ref{eqn:deffp}),(\ref{eqn:definitionofg}) and (\ref{eqn:simplificationofgT}) put the main identity (\ref{eqn:maintheorem}) in
theorem \ref{thm:main} as a sum over all homotopy classes of essential
embedded 3-holed spheres in the surface $F$.  At this moment, we
are not able to reconcile the two different expressions in (\ref{eqn:deffp})
and (\ref{eqn:simplificationofgT}). The function $\L(\frac{1}{\cosh^2(x/2)})$ was first
introduced and used by Bridgeman \cite{bridge}.

%There is a high possibility that the above functions $f,g$ can be
%simplified drastically using the pentagon relation for dilogarithm
%function. We are currently working on it.

%The identity (\ref{eqn:maintheorem}) in Theorem \ref{thm:main} now holds with $f$ and $g$ defined as above.

%\bigskip

For the identity (\ref{eqn:maintheorem2}) in Theorem
\ref{thm:mainwithboundary}, the functions $\hat f$ and $\bar f$
are defined using the lasso function $La(l,m)$ and the function
$f(P)$  as follows:

\begin{eqnarray*}\label{eqn:definitionofhatfandbarf}
  %f(P) &=& 4\pi^2 -\big( 8\sum_{i=1}^3\big( b(m_i)+ b(p_i)\big)+4\sum_{i\neq j}La(l_i,m_j)\big) \\
  {\hat f}(P) &:=& f(P)+8\big(\L(\frac{1}{\cosh^2(p_1/2)})+La(l_2,m_3)+La(l_3,m_2)\big)
\end{eqnarray*}
where $\partial P \cap \partial F =L_1$ and
\begin{eqnarray*}
  {\bar f}(P) &:=& f(P)+8(\L(\frac{1}{\cosh^2(p_1/2)})+ \L(\frac{1}{\cosh^2(p_2/2)})+ \L(\frac{1}{\cosh^2(m_3/2)})) \\
              &+& 8\big(La(l_2,m_3)+La(l_3,m_2)+La(l_3,m_1)+La(l_1,m_3) \big)
\end{eqnarray*}
where  $\partial P \cap \partial F=L_1 \cup L_2$.

\bigskip

\noindent {\bf Remark.} The expressions $f(P)$, $g(T)$, $\hat
f(P)$ and $\bar f(P)$ defined above are still valid if $P$ or $T$
are hyperbolic surfaces with some cusp ends. Namely, if some $l_i$
or $c$  tend to 0 (which imply the corresponding $m_j$'s and
$p_i$'s tend to infinity), the functions $f,g,\bar f, \hat f$
converge to well defined limit functions.
%are zero in the above formulas, and
%some of the $m_i$'s are infinite, the functions are still well
%defined, as the limit is well defined, when some of the boundary
%lengths approach 0.
If we use these limit functions in (\ref{eqn:maintheorem2}), then (\ref{eqn:maintheorem2}) becomes the
identity for finite area hyperbolic surfaces $F$ with geodesic
boundary and cusp ends. In this case, the right-hand-side of (\ref{eqn:maintheorem2})
is the volume of the unit tangent bundle of $F$ and the
left-hand-side is the sum over all hyperbolic 3-holed spheres $P$
and 1-holed torus $T$ where $P$ may have cusp ends. For
simplicity, we omit the details here. Some details, including an
identity for the cusped torus can be found in \cite{HuTan}.

%{\bf Feng: could you add the similar formulas for $f_i$,
%$f_{i,j}$? Thanks}
%{\bf Tan: 1. Ok, done.
%  2.I've also made some corrections to the formula for $f$, for example $p_i$ instead of $b_i$ and some of the brackets were also readjusted.}

%% file: sec3.tex
\section{Decomposing the unit tangent bundle of the surface}\label{s:decomposing}
%{\bf Tan's comments: We should adopt a uniform notation for the unit tangent bundle of a surface, I notice in introduction you used $S(F)$, which is good except I used $S$ for simple surfaces, anyhow, it is fine with me whichever we use, here it's still what I used which is $S(T_{F})$.}

Suppose $F$ is a compact hyperbolic surface with or without
boundary so that if $\partial F \neq \emptyset$, then $\partial F$
consists of geodesics. Let $S(F)$ be the unit tangent bundle of
$F$ and  $\mu$ be the measure on $S(F)$ invariant under the
geodesic flow so that $\mu(S(F))=-4\pi^2 \chi(F)$.

We will produce a decomposition of $S(F)$ as follows. Given a
vector $v \in S(F)$, let $g_v^+$ and $g_v^-$ be the geodesic rays
determined by $v$ and $-v$. By the ergodicity of the geodesic
flow, for generic choice of $v$ with respect to $\mu$, we may
assume that

(1) if $\partial F =\emptyset$, each geodesic ray $g^-_v$ and
$g_v^+$ is not simple and intersects every closed geodesic,

(2) if $\partial F  \neq \emptyset$,  each geodesic ray $g^-_v$
and $g_v^+$ intersects $\partial F$.

Indeed, the set $X$ of all $v$'s in $S(F)$ satisfying (1) is
invariant under the geodesic flow. Furthermore, the set $X$ has
positive $\mu$-measure. It follows that $\mu(S(F) -X)=0$.  To see
(2), we apply the ergodicity of the geodesic flow to the metric
double of $F$ across the boundary of $F$.

In the sequel, we will focus only on these generic vectors $v$.

Given a generic vector $v \in S(F)$, we define an associated graph
$G(v)$ to $v$ as follows.

Let $t_1>0$ be the smallest number so that the geodesic segment
$g_v^-[0,t_1] \cup g_v^+([0, t_1])$ either intersects $\partial F$
or intersects itself. Say for simplicity that this occurs in the
ray $g^+_v$. This means $g_v^-([0, t_1]) \cup g_v^+([0, t_1))$ is
a simple path in $F$ so that $g_v^+(t_1)$ is in $\partial F$ or in
$g_v^-([0, t_1]) \cup g_v^+([0, t_1))$.  Next, let $t_2 \geq t_1$
be the smallest number so that $g_v^-(t_2)$ is either in $\partial
F$ or in $g_v^-([0, t_2)) \cup g_v^+([0, t_1])$.  The associated
graph $G(v)$ is defined to be the connected component of
$g_v^-([0, t_2]) \cup g_v^+([0, t_1])\cup
\partial F$ which contains $g_v^+(0)$.  In particular, $G(v)
=g_v^-([0, t_2]) \cup g_v^+([0, t_1])$ if $\partial F =\emptyset$.
By the construction, the Euler characteristic of $G(v)$ is always
$-1$. We also define $G(v)^o:=g_v^-([0, t_2]) \cup g_v^+([0, t_1])$, with the
orientation induced from $v$. Note that $G(v)^o=G(v)$ if $G(v) \cap \partial F =\emptyset$, otherwise, it
 is a strict subset of $G(v)$. See figure 1 for closed surfaces and
figure 4 for surfaces with non-empty boundary.

\begin{figure}[ht!]
\centering
\includegraphics[scale=0.55]{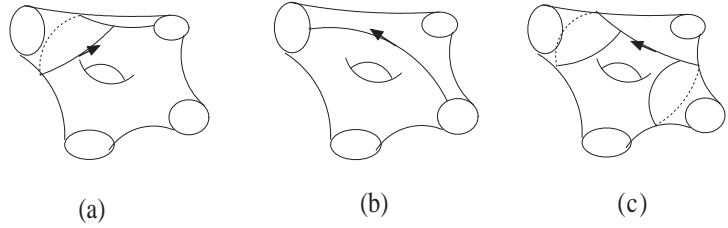}
\caption{creation of spine} \label{figure 4}
\end{figure}

Recall that for a hyperbolic surface $F$,  a compact embedded
subsurface $S \subset F$ is said to be \emph{ geometric }if the
boundaries of $S$ are geodesic and \emph{proper} if the inclusion
map $i: S \rightarrow F$ is injective. Furthermore a surface is
simple if it is a 3-holed sphere or 1-holed torus.

\begin{prop}\label{prop:loopisembedded}
The graph $G(v)$ is contained in a unique geometric embedded
simple surface $\Sigma(v)$.
%Furthermore, $G(v)$ is a deformation retract of $\Sigma(v)$.
\end{prop}

\begin{proof}
%\noindent{\it Proof}. %We consider the case where $\partial F =\emptyset$.
%({\bf Feng:  I am not so sure about the other cases mainly because
% for $\partial F \neq \emptyset$, the sets $\gamma_v$ and $\gamma_v^{tr}$ have to be defined
% differently, taking care of the boundary. See my remarks in \S4.1. We would drop
% the sentence "other case is similar"}).
% {\bf Tan: Agreed, in fact, given the previous paragraph, I am removing it.}
Cutting $F$ open along $G(v)$, we obtain a (possibly disconnected) surface whose metric
completion $\hat F$ is a (possibly disconnected) compact hyperbolic surface with convex
boundary. The boundary of $\hat F$ consists of simple closed
geodesics (corresponding to components of $\partial F$ not in $G(v)$) and piecewise simple
geodesic loops (corresponding to $G(v)$).

If  $\hat \gamma$ is a piecewise simple geodesic loop in $\partial
\hat F$, it is freely homotopic to a simple closed geodesic
$\gamma$ in $\hat F$ which is a component of the boundary of  the
convex core $core( \hat F)$ of $ \hat F$.  Furthermore $\hat
\gamma$ and $\gamma$ are disjoint by convexity. Therefore, $\hat
\gamma$ and $\gamma$ bound a convex annulus exterior to $core(\hat
F)$ and $G(v)$ is disjoint from  $core(\hat F)$.
The subsurface $\Sigma(v)\subset F$ is the union of these convex
annuli bounded by $\hat \gamma$ and $\gamma$. The Euler
characteristic of $\Sigma(v)$ is $-1$ by the construction.
The surface $\Sigma(v)$ is
 unique. Indeed, if
 $\Sigma' \neq \Sigma \subset F$ is a simple geometric subsurface
 so that $G(v) \subset \Sigma'$, then $\Sigma'$ has a
 boundary component say $B$ which
 intersects one of the boundaries $ \gamma$ of $\Sigma$ transversely.
 Therefore, $B$ must intersect the other boundary
 $\hat \gamma$ of the convex annulus described earlier. Hence
 it intersects $G(v)$ which contradicts $G(v)\subset \Sigma'$.
  \end{proof}
% The
%shortest distance retraction of $\hat F$ to $core(F^*)$ produces a
%retraction from $\Sigma(v)$ to $G(v)$. \qed

Note that topologically a regular neighborhood  $N(G(v))$ of the
graph $G(v)$ is either the 3-holed sphere $F_{0,3}$ or the 1-holed
torus $F_{1,1}$.  In the case that $N(G(v)) \cong F_{0,3}$ so that
two boundary components of $N(G(v))$ are freely homotopic, then
$\Sigma(v) \cong F_{1,1}$ and $\Sigma(v) -G(v)$ consists of two
annuli (see figure 1(c)). In this case, $G(v)$ is not a
deformation retract of $\Sigma(v)$. In all other cases,
$\Sigma(v)$ is isotopic to $N(G(v))$ so that $\Sigma(v)$
deformation retracts to $G(v)$.

 As a consequence,
we have produced the following decomposition of the unit tangent
bundle $S(F)$. Given a simple geometric subsurface $\Sigma$ in
$F$, let

$$W(\Sigma) = \{ v \in S(F) | G(v) \subset \Sigma \}.$$

Then by proposition \ref{prop:loopisembedded}, we have the following decomposition
$$
S(F) = \bf{Z} \bigsqcup \bigsqcup_{P} W(P) \bigsqcup_{T} W(T)
$$
where $\bf{Z}$ is a set of measure zero and the union is over all
simple geometric 3-holed spheres $P$ and  1-holed tori $T$.

Take the $\mu$ measure of the above decomposition, we obtain the
main identities (\ref{eqn:maintheorem}) and (\ref{eqn:maintheorem2}) in Theorems \ref{thm:main} and \ref{thm:mainwithboundary}
\begin{equation}
\mu(S(F)) = \sum_{P}\mu(W(P)) + \sum_{T} \mu(W(T)).
\end{equation}
The focus of the rest of the paper is to calculate the volume of
$W(\Sigma)$ for simple surfaces $\Sigma$.

We end this section with a related simpler decomposition of $S(F)$
indexed by the set of all simple closed geodesics.  For
simplicity, we assume that $F$ is a closed hyperbolic surface.
Given a generic unit tangent vector $v$, the geodesic ray $g_v^+$
intersects itself. Let $t_1>0$ be the first time so that
$g_v^+(t_1) \in g_v^+([0, t_1))$, say $g_v^+(t_1) = g_v^+(t_2)$
for some $0 \leq t_2 < t_1$. Then $g_v^+|_{[t_2, t_1]}$ is a
simple loop freely homotopic to a simple closed geodesic $s$ in
$F$. Denote $g_v^+([t_2, t_1])$  by $Lop(v)$. For any given simple
closed geodesic $s$ in $F$, let $U(s) =\{ v \in S(F)| Lop(v) \cong
s$\}. Then we obtain a decomposition $S(F) = \mathbf{Z'} \bigsqcup
\bigsqcup_s  U(s)$ where the disjoint union is indexed by the
simple closed geodesics $s$ and $\mu(\mathbf{Z'})=0$. The
associated identity is $\mu(S(F)) = \sum_{s} \mu(U(s))$. However,
we are not able to calculate $\mu(U(s))$. It is not clear if
$\mu(U(s))$ depends only on the length of $s$ and the topology of
$F$.

%% file: sec4.tex
\section{Identifying the sets in the decomposition}\label{s:measureofdecompostion}

We will investigate the sets  $W(P)$ and $W(T)$ by studying their
complements in $S(P)$ and $S(T)$. We will decompose the
complementary sets into a disjoint union of sets, and identify
each with subsets of $S(\HH)$ in this section so that the
computation of their volume can be
 carried out in \S5.

 For simplicity, we will deal with closed hyperbolic surfaces
 $F$. The modification for surfaces with non-empty boundary is
 easy, see \S\ref{ss:surfaceswithboundary}.  For a generic unit tangent vector $v \in S(F)$, $G(v)$
 is a graph lying in a simple geometric surface $\Sigma$ of $F$. In
 particular, $v \in S(\Sigma)$. Now for $v \in S(\Sigma)$, let $G_{\Sigma}(v)$ be the associated
 graph of $v$ in $\Sigma$. Then by definition,
 $$ W(\Sigma) = \{ v \in S(\Sigma) | G(v) = G_{\Sigma}(v)\}.$$
 To calculate $\mu(W(\Sigma))$, we will focus on the complement
 $$ V(\Sigma) = S(\Sigma) -W(\Sigma)=\{ v \in S(\Sigma) | G_{\Sigma}(v) \cap \partial \Sigma \neq \emptyset\}
 .$$

To this end, recall that $G_{\Sigma}(v) = g_v^+([0,t_1]) \cup
g_v^-([0,t_2]) \cup B$ where $B$ consists of some boundary
components of $\Sigma$ ($B$ could be the empty set). Let
$G_{\Sigma}^+(v)$ and $G_{\Sigma}^-(v)$ be the geodesic paths
$g_v^+|_{[0, t_1]}$ and $g_v^-|_{[0, t_2]}$ defined in \S\ref{s:decomposing}. By
definition,
$$ V(\Sigma) = \{ v \in S(\Sigma) | \quad G_{\Sigma}^+(v) \quad \text{or}  \quad
 G_{\Sigma}^-(v) \quad \text{is a simple arc ending at} \quad \partial \Sigma\}.$$

There are two cases which can occur. Namely either both
$G_{\Sigma}^+(v)$ and $G_{\Sigma}^-(v)$ are simple arcs ending at
$\partial \Sigma$ or exactly one of them ends at $\partial
\Sigma$.

These two cases will be discussed separately in the case $\Sigma$ is a 3-holed sphere
in the subsections
\S\ref{ss:bothsimplearcs} and \S\ref{ss:lassos} below, and in \S\ref{ss:vinT} in the case $\Sigma$ is a 1-holed torus. We will first recall some facts in \S\ref{ss:prelimonconvexsurfaces}
about convex hyperbolic surfaces.

\subsection{Preliminaries on convex surfaces}\label{ss:prelimonconvexsurfaces}

Suppose $X$ is a compact connected surface with a hyperbolic
metric so that $\partial X $ consists of convex curves. Then,
unless $X$ is simply connected, each component of
$\partial X$ is an essential loop in $X$ homotopic to a geodesic.
As a convention, we will identify the universal cover $\tilde{X}$
of $X$ with a convex subset of $\HH$.

The following notation and conventions will be used. For a
hyperbolic surface $Y$, a \it geodesic path \rm is a map $s: [a,b]
\to Y$ satisfying the geodesic equation so that $s'(t) \in S(Y)$.
We are mainly interested in geodesic paths whose end points are in
$\partial Y$.  A geodesic path $s$ is called a \it geodesic loop
\rm if $s(a)=s(b)$. A \it simple \rm geodesic path or loop
satisfies the condition that $s|_{(a,b)}$ is an injective map. Two
paths $\alpha_i: ([a_i, b_i], \{a_i, b_i\}) \to (X,
\partial X)$, $i=0,1$, are \it homotopic, \rm denoted by $\alpha_0
\cong \alpha_1$ if there is a homotopy $H: ([0,1] \times [0,1],
\{0,1\}\times [0,1]) \to (X,
\partial X)$ so that $H(t, i) = \alpha_i(a_i + t(b_i-a_i))$ for
$i=0,1$ and all $t$.  Two loops $\alpha_i$, $i=0,1,$ with the same
 base point $p=\alpha_i(a_i)=\alpha_i(b_i)$ which are \it
relatively homotopic with respect to $p$ \rm will be denoted by
 $\alpha_0 \cong \alpha_1$ rel\{p\}.

The main technical result in this subsection is the following:

\begin{prop}\label{prop:basicfactssimplegeodesics}
Suppose $X$ is a compact non-simply connected hyperbolic surface
with convex boundary.
\begin{enumerate}
  \item If $X$ is a topological annulus, then any geodesic path $s$ in
  $X$ joining different boundary components of $X$ is simple;
  \item  If $s \cong t$ are two geodesic paths in $X$ joining
  different boundary components of $X$ and $t$ is simple, then $s$
  is simple;
  \item If $p \in \partial X$ and $s:([0, a], \{0, a\}) \to (X, \{p\})$  is a geodesic
  path so that $s \cong t $ rel\{p\}  and $t$ is a simple loop, then $s$
  is a simple loop.
\end{enumerate}
\end{prop}

\begin{proof}
We will need the following simple lemma whose proof is omitted.

\begin{lem}\label{lem:nonintersectinggeodesics}  Suppose $\gamma$ is a hyperbolic isometry of
$\HH$ with axis $A$ and $g$ is a geodesic intersecting $A$
transversely.   Then $\gamma^n (g) \cap g =\emptyset$ for all $n
\in \mathbf{Z}-\{0\}$.
\end{lem}

%Indeed, after a conjugation, we may assume that in the
%upper-half-plane model of the hyperbolic plane, $A$ is the y-axis
%and $\gamma(z) = kz$ for some $k >1$. As the geodesic $g$
%intersects $A$ transversely, it has end points $\{a,b\}$ where
%$a<0$ and $b>0$. It follows that $\gamma^n(g)$ has end points
%$k^na$, $k^nb$ and therefore $\gamma^n (g) \cap g =\emptyset$.

To see (1), let $c$ be the unique simple closed geodesic in $X$.
Then $s$ must intersect $c$ in $X$. Lifting $s$ and $c$ to the
universal cover and using the above lemma, we see that any two
distinct lifts of $s$ in $\tilde{X}$ are disjoint. Thus $s$ is
simple.

To see (2), suppose otherwise there exist two distinct lifts $s_1$
and $s_2$ of $s: [0, d] \to X$  in $\tilde{X}$ so that the
interiors of $s_1$ and $s_2$ intersect.  Let $s$ join boundary
components $a$ and $b$ of $X$ and $\tilde{a_i}$ and $\tilde{b_i}$
be the lifts of $a$ and $b$ so that $s_i(0) \in \tilde{a_i}$ and
$s_i(d) \in \tilde{b_i}$. Since $s \cong t$, by the homotopy
lifting theorem, there exist two distinct lifts $t_1$ and $t_2$ of
$t$ in $\tilde{X}$ so that $t_i$ joins $\tilde{a_i}$ to
$\tilde{b_i}$.

\begin{figure}[ht!]
\centering
\includegraphics[scale=0.5]{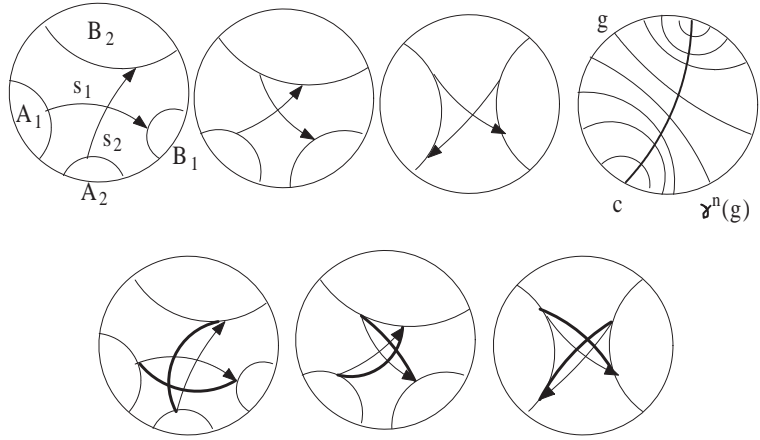}
\caption{lifting and disjointness} \label{figure 5}
\end{figure}

We claim that interiors of $t_1$ and $t_2$ intersect. This in turn
contradicts the fact that $t$ is simple.

To see the claim, first, we note that  $\tilde{a_1}$ is disjoint
from $\tilde{a_2}$. For otherwise, $s_2 =\gamma^n (s_1)$ for a
deck transformation element $\gamma$ corresponding to the boundary
$a$ of $X$. Furthermore, due to convexity both $s_1$ and $s_2$
intersect the axis of $\gamma$.  Thus  by the lemma above, $s_1$
is disjoint from $s_2$ which contradicts the assumption. By the
same argument we see that $\tilde{b_1}$ is disjoint from
$\tilde{b_2}$. Since $a \cap b =\emptyset$ by assumption, we see
that $\{\tilde{a_1}, \tilde{a_2}, \tilde{b_1}, \tilde{b_2}\}$
consists of four distinct convex curves  in $\tilde{X}$. Let $A_1,
A_2, B_1, B_2$ be the four disjoint half-spaces in $\HH$ bounded
by these four convex curves. Let $S^1_{\infty}$ be the circle at
infinity of the hyperbolic plane. Then $s_1 \cap s_2 \neq
\emptyset$ is equivalent to saying that $\overline{A_1} \cap
S^1_{\infty}$ and $\overline{B_1} \cap S^1_{\infty}$ are in the
different components of $S^1_{\infty} -\overline{A_2} \cup
\overline{B_2}$. This in turn implies that interiors  $t_1$ and
$t_2$ intersect. Thus part (2) holds.

The proof of part (3) is similar to that of (2). Suppose the
result is false. Then there exist two distinct lifts $s_1, s_2:
[0, d] \to \tilde{X}$ of $s: ([0, d], \{0, d\}) \to (X, \{p\})$ so
that $s_1(d_1) = s_2(d_2)$ for some $d_1, d_2 \in (0, d)$. Let
$t_1, t_2$ be the lifts of $t$ so that the end points of $t_i$ are
the same as that of $s_i$ (by the homotopy lifting property). We
claim that the interior of $t_1$ intersects the interior of $t_2$.
This would produce a contradiction to the assumption on $t$.

To see this, let $s_i(0) \in \tilde a_i$ and $s_i(d) \in \tilde
b_i$ where $\tilde{a_i}$ and $\tilde{b_i}$  are lifts of the same
boundary $a$ of $X$. By the same argument as above, we see that
$\tilde{a_1} \cap \tilde{a_2} =\emptyset$ and $\tilde{b_1} \cap
\tilde{b_2} =\emptyset$. However, it is possible that $\tilde a_1 =
\tilde b_2$ and $\tilde a_2=\tilde b_1$. If $\{\tilde{a_1},
\tilde{a_2}, \tilde{b_1}, \tilde{b_2}\}$ are pairwise disjoint,
then the same argument as above shows that the claim holds. In the
other cases, $\{\tilde{a_1}, \tilde{a_2}, \tilde{b_1},
\tilde{b_2}\}$ consists of 2 or 3 geodesics, the same argument
again shows that the interiors of $t_1$ and $t_2$ intersect since
$s_1$ and $s_2$ are geodesics and $t_i$ and $s_i$ have the same
end points. See figure \ref{figure 5}.
\end{proof}

\bigskip
\subsection{ Vectors $v$ in $V(P)$ so that $G_{P}^+(v)$
and $G_{P}^-(v)$ are simple arcs ending at $\partial P$}\label{ss:bothsimplearcs}

We begin by recalling the beautiful work of M. Bridgeman
\cite{bridge} relevant to our setting. Given a compact hyperbolic
surface $X$ with geodesic boundary and a (not necessarily simple)
geodesic path $\alpha: ([0, a], \{0, a\}) \to (X,
\partial X)$  so that $\alpha'(0)$ and $\alpha'(a)$ are
perpendicular to $\partial X$, let
$$  H(\alpha)  =\{ s'(t)  \\ | s: ([0, b],  \{0,b\}) \to (X,
\partial X) \\ \text{ geodesic, so that } s \cong \alpha\}.$$

\begin{thm}(Bridgeman)
The measure $\mu(H(\alpha))$ of $H(\alpha)$ is
$4\L(\frac{1}{\cosh^2(l(\alpha)/2)})$ where $l(\alpha)$ is the
length of $\alpha$.
\end{thm}

Calegari gave a very nice short and elegant proof of this in \cite{cal}. If
we use $\alpha^{-1}$ to denote the reversed path $\alpha^{-1}(t) =
\alpha( a-t)$, then the measures of $H(\alpha^{-1})$ and
$H(\alpha)$ are the same. In Bridgeman's work, he considered
unoriented paths, i.e., the elements in $H(\alpha) \cup
H(\alpha^{-1})$ and showed that its measure is
$8\L(\frac{1}{\cosh^2(l(\alpha)/2)})$.  For simplicity, we use
$H(\alpha^{\pm 1})$ to denote $H(\alpha) \cup H(\alpha^{-1})$.

The main result in this section is to prove:

\begin{prop}\label{prop:bridgemansetsinP} Suppose $P$ is a hyperbolic 3-holed 3-sphere with geodesic
boundary components $L_1, L_2, L_3$ and the shortest paths joining
boundary components being $M_i$ and $B_i$ as in \S\ref{ss:pop}. Then
\begin{enumerate}
  \item $ \{ v \in S(P) | G_{P}^+(v), G_{P}^-(v) \text{ both simple arcs ending at
  $\partial P$} \}
   \subset  \cup_{ i=1}^3 ( H(M_i^{\pm 1})
     \cup H(B_i^{\pm 1})).$
  \item
$\cup_{ i=1}^3 ( H(M_i^{\pm 1}) \cup
  H(B_i^{\pm 1})) \subset V(P).$
\end{enumerate}
\end{prop}

\begin{proof}  To see (1), by the construction of $G_{P}^+(v)$ and
$G_{P}^-(v)$, the interiors of these two simple arcs are disjoint.
It follows that the geodesic path $G_{P}^+(v) \cup G_{P}^-(v)$ is
a simple path with end points in $\partial P$.
%In the case that $G_{\Sigma}^+(v) \cup G_{\Sigma}^-(v)$ is a
%simple loop, we can produce a homotopy relative to $\partial P$ so
%that $G_{\Sigma}^+(v) \cup G_{\Sigma}^-(v)$ is homotopic to a simple arc.
It is well known that any simple path $s: ([0, 1], \{0, 1\}) \to
(P,\partial P)$ is homotopic to  $M_i$, or $B_i$, or a point. The
 path $ G_{P}^+(v) \cup G_{P}^-(v)$ cannot be homotopic to
a point since it is a geodesic path. Thus the conclusion follows.

To see (2),  let $s:[0, a] \to P$ be a geodesic path homotopic to
$M_i$ or $B_i$. If $s \cong M_i$, by proposition
\ref{prop:basicfactssimplegeodesics} (2), $s$ is simple. Thus
$s'(t) \in V(P)$. If $s \cong B_i$, we claim that there exists $b
\in (0, a)$ so that $s|_{[0,b]}$ and $s|_{[b,a]}$ are simple arcs.
To see this, first of all, the path $s$ intersects $M_i$ in
exactly one point. Indeed, if there are at least two points of
intersection, then there will be a lift $\tilde{s}$ of $s$ in the
universal cover $\tilde P$ so that $\tilde{s}$ intersects two
distinct lifts $a_1$ and $a_2$ of $M_i$. Let $\tilde{B}$ be the
lift of $B_i$ so that both $\tilde{B}$ and $\tilde{s}$ start and
end at the same geodesics  which are lifts of $L_i$. Then
$\tilde{B}$ intersects $a_1$ and $a_2$, i.e., $B_i$ intersects
$M_i$ at two points. This is impossible. Furthermore, by
topological reasons, $s$ must intersect $M_i$. It follows that $s$
intersects $M_i$ in exactly  one point, say $s(b) \in M_i$ for
some $b \in (0, a)$. We claim that $\alpha:=s|_{[0,b]}$ and $\beta:=s|_{[b,a]}$
are both simple arcs. Indeed, let $X$ be the surface obtained by
cutting $P$ open along $M_i$. Then $X$ is a convex hyperbolic
surface homeomorphic to an annulus. Both paths $s|_{[0,b]}$ and
$s|_{[b,a]}$ are geodesics in $X$ joining different boundary
components of $X$. Thus, by proposition \ref{prop:basicfactssimplegeodesics}, both of them are
simple and both intersect the unique closed geodesic $L_i$ in $X$.
It is well known that if $x,y$ are two oriented geodesics in a
convex hyperbolic annulus $X$ so that both $x,y$ intersect the
closed geodesic in $X$,  then all intersection points between $x$
with $y$ have the same intersection sign. See figure \ref{figure 6}(a) for a
pictorial explanation in the universal cover.  Thus all
intersection points between $\alpha$ and $\beta$ have the same
sign.

%Furthermore, we note that all intersection points between
%$\alpha=s|_{[0,b]}$ and $\beta=s|_{[b,a]}$ have the same algebraic
%intersection sign counted from $\alpha$ to $\beta$. Indeed
%$\alpha$ and $\beta$ are oriented geodesics in the complete
%hyperbolic cylinder $X^*$ which contains $X$ as a convex
%subsurface and $\alpha$, $\beta$ intersect the unique closed
%geodesic $a^*$ in $X^*$. Now it is well known that if $x,y$ are
%two oriented geodesics in $X^*$ so that both $x,y$ intersect
%$a^*$, then all intersection points between $x$ with $y$ have the
%same intersection sign. See figure 6(a) for a pictorial
%explanation in the universal cover.  Thus all intersection points
%between $\alpha$ and $\beta$ have the same sign.

We now finish the proof of (2) by showing that for any $t \in (0,
a)$, $s'(t) \in V(P)$.  Suppose otherwise, there exists $t_0\in (0,
a)$ so that
$v = s'(t_0) \in W(P)$, i.e., the graph $G_P(v)$ does not
intersect $\partial P$. Since $s([0, a])$ is a union of two simple
arcs, the graph $G_P(v)$, considered as a sub-path $s|_{[T_1,
T_2]}$ in $s$, is a union of two simple arcs. Since $G_P(v)$ is
embedded in the planar surface $P$, by definition of $G_P(v)$,
there are two possible embedding of $G_P(v)$ in $P$ as shown the
figure \ref{figure 6}(b),(c).

\begin{figure}[ht!]
\centering
\includegraphics[scale=0.55]{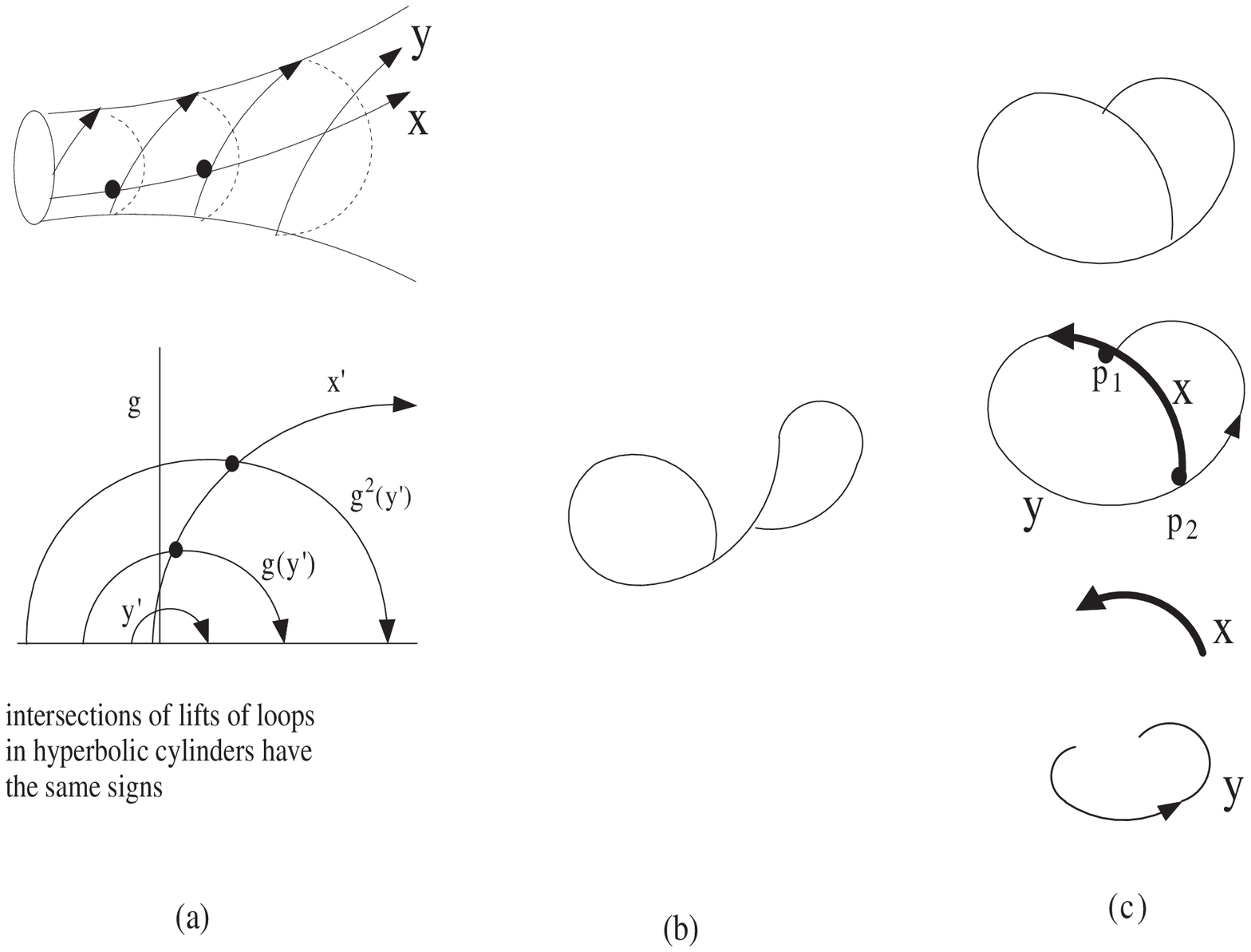}
\caption{embedding of graphs} \label{figure 6}
\end{figure}

In the first case, there are two disjoint simple loops in the
graph $G_P(v)$. In this case, $G_P(v)$ cannot be a union of two
simple arcs due to the disjoint simple loops.  In the second case,
the graph $G_P(v)$ can be expressed as a union of two simple arcs
$x$ and $y$ in an essentially  unique way as shown in figure \ref{figure 6}(c).
Let $p_1$ and $p_2$ be the two vertices of $G_P(v)$ and orient
both arcs $x$ and $y$. Then the intersection signs at $p_1$ and
$p_2$ from $x$ to $y$ are opposite. It follows that this case does
not occur in $s([0,a])$ by the calculation above. This ends the
proof of (2).

\end{proof}

\subsection{ Lassos}\label{ss:lassos}

For a hyperbolic 3-holed sphere $P$, it remains to identify
 the set $V(P)-\bigcup_{i=1}^3 (H(M_i^{\pm 1}) \cup H(B_i^{\pm 1}))$. If
 $v$ is in the set, then one of $G_P^+(v)$ or $G_P^-(v)$ is a
 simple arc ending at $\partial P$ and the other one is part of, or contains
  a loop. Thus $G_P(v)$ is a lasso (see figure \ref{figure 4}(a)), as defined below.

\begin{defn} (Lassos)\label{def:lassos}
Let $X$ be a hyperbolic surface with geodesic boundary. A
positively oriented \emph{lasso} on $X$ is a geodesic path
$$\alpha:[T_1,T_2] \rightarrow X$$ such that
\begin{enumerate}
  \item  $\alpha(T_1) \in \partial X$,
  \item  $\alpha$ is injective on $(T_1,T_2)$, and
  \item  $\alpha(T_3)=\alpha(T_2)$ for some $T_1 \leq T_3 < T_2$.
\end{enumerate}
The image of $\alpha$, ignoring orientation, is a lasso.
A negatively oriented lasso $\beta$ is a geodesic path so that
$\beta(-t)$ is a positively oriented lasso.  Call $\alpha(T_1)$ the
\emph{base point}, $\alpha(T_2)=\alpha(T_3)$ the \emph{knot},
$\alpha[T_1, T_3]$ the \emph{stem}, $\alpha|_{[T_3, T_2]}$ the
\emph{loop}, and
 $\alpha(\frac{T_2+T_3}{2})$ the \it midpoint \rm of
the loop of the lasso. Note that $\alpha(T_1, T_2) \cap \partial X =
\emptyset$.

\end{defn}

The midpoint of the loop $\alpha(\frac{T_2+T_3}{2})$ is diametrically opposite
to the knot in the loop of a lasso.  If $\gamma$ is the unique
oriented geodesic on $X$ homotopic to the loop of $\alpha$, then
%by cutting along $\alpha$ and using a convexity argument as before,
% they are disjoint, and
the loop of $\alpha$ and $\gamma$ bound a hyperbolic cylinder $A$
embedded in $X$. It is easy to see by lifting to the universal
cover that $\alpha(\frac{T_3+T_2}{2})$ is the point on the loop
which is closest to $\gamma$ on the cylinder $A$.

Note that if $\alpha$ and $\beta$ are two lassos so that $\alpha$
is positively oriented and $\beta$ is negatively oriented, then by
definition $\alpha'(t) \neq \beta'(t')$ for all parameters $t,
t'$.  Furthermore, the involution map $v \to -v$ in $S(X)$ sends
tangent vectors to positively oriented lassos to that of
negatively oriented lassos.  Thus it suffices to calculate the
measure of tangents to positively oriented lassos.

\begin{prop}\label{prop:vonlasso}
Suppose that $\alpha:[T_1,T_2] \rightarrow \Sigma$ is a positively
oriented lasso in a compact hyperbolic surface $\Sigma$ with
geodesic boundary, and $\alpha(T_3)=\alpha(T_2)$ is the knot of
$\alpha$. Then $G_{\Sigma}(\alpha'(t))- \partial
\Sigma=\alpha((T_1, T_2])$  if and only if
$T_1 \le t \leq \frac{T_2+T_3}{2}$.
\end{prop}
\begin{proof} The midpoint  $\alpha(\frac{T_2+T_3}{2})$ lies on the critical set where if we
exponentiate in both directions at equal speed, we reach the knot
of the lasso at the same time. For $T_1<t<\frac{T_2+T_3}{2}$, we
get the lasso $\alpha$ and for $\frac{T_2+T_3}{2}<t<T$, we will
exponentiate in the other direction of the knot and
$G_{\Sigma}(\alpha'(t))$ will not include the stem of $\alpha$.
\end{proof}

\bigskip

In the rest of the discussion, we assume surfaces are oriented so
that their boundaries have the induced orientation.
 Given $v \in V(P)
-\bigcup_{i=1}^3 (H(M_i^{\pm 1}) \cup H(B_i^{\pm 1}))$, the graph
$G_P(v)^o$ (see \S\ref{s:decomposing}) is a lasso. Since its loop is simple, it is freely
homotopic to $L_i^{\pm 1}$ for some $i$.

For $i,j,k$ distinct, let $W(L_i, M_j)$  be $\{ v \in S(P)|$
$G_P(v)^o$ is a positive lasso whose loop is homotopic to $L_i$,
the base point of $G_P(v)$ is in $L_k$, and $v \notin
\bigcup_{l=1}^3 (H(M_l^{\pm 1}) \cup H(B^{\pm 1}_l))\}$. Let
$W(L_i^{-1}, M_j)$ be the set defined in the same way except the
loop of the lasso is homotopic to $L_i^{-1}$. Let $\mathbf{A}:
S(P) \to S(P)$ be the involution $\mathbf{A}(v) = -v$.
$\mathbf{A}$ sends vectors generating positive lassos to those
generating negative lassos, and vice versa, since $G_P(v)^o=G_P(-v)^o$
with opposite orientations.

We have:

\begin{lem}\label{lem:VPdecomposition}
  The set $V(P) -\bigcup_{i=1}^3
(H(M_i^{\pm 1}) \cup H(B_i^{\pm 1}))$ can be decomposed as
$$\bigsqcup_{ i \neq j} (W(L_i, M_j) \cup W(L_i^{-1}, M_j))
\bigsqcup \mathbf{A}(\bigsqcup_{ i \neq j} (W(L_i, M_j) \cup
W(L_i^{-1}, M_j)) ).$$ In particular,  $$\mu(V(P)) = 8
\sum_{i=1}^3 (\L(\frac{1}{\cosh^2(m_i/2)}) +
\L(\frac{1}{\cosh^2(p_i/2)})) + 4 \sum_{ i \neq j} \mu(W(L_i,
M_j)).$$
\end{lem}

\begin{proof} The decomposition in the first sentence follows from the above discussion. We claim that $W(L_i, M_j)$ and $W(L_i^{-1}, M_j)$ are related by
an isometry of $P$. Indeed, the hyperbolic 3-holed sphere $P$
admits an orientation reversing isometry $R$ so that $R|_{M_i} =
id$ and $R$ interchanges the two hexagons obtained by cutting $P$
open along $M_i$'s.  In particular, $R$ reverses the orientation
of each boundary component. Therefore, the derivative $R_*$ of $R$
sends $W(L_i, M_j))$ to $W(L_i^{-1}, M_j)$, i.e., $R_*(W(L_i,
M_j)) =W(L_i^{-1}, M_j))$. In particular, $\mu(W(L_i, M_j)) =
\mu(W(L_i^{-1}, M_j))$.
\end{proof}

\bigskip
\subsection{ Understanding the set $W(L_i, M_j)$}\label{ss:understandingWLM}
\medskip
We begin with some notation. The circle at infinity of the
hyperbolic plane is denoted by $S^1_{\infty}$. Given $x \neq y \in
\HH \cup S^1_{\infty} $, let $G[y,x]$ be the oriented geodesic
from $y$ to $x$.  In particular, if $x \neq y \in S^1_{\infty} $,
then $G[y,x]$ is the complete oriented geodesic determined by $y,x$.

\begin{figure}[ht!]
\centering
\includegraphics[scale=0.7]{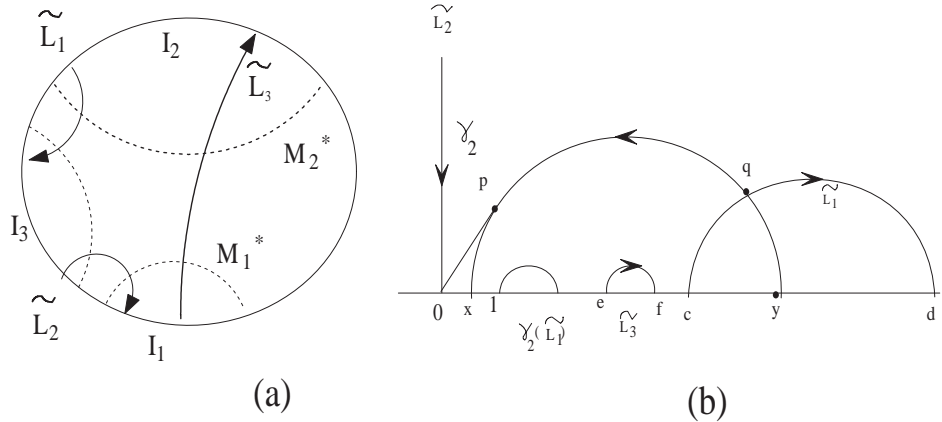}
\caption{lifts of boundary} \label{figure 7}
\end{figure}

Consider the universal cover $\tilde P$ of the hyperbolic 3-holed
sphere $P$ as a convex subset of $\HH$ so that the covering map is
$\Pi: \tilde P \to P$. We assume that $\tilde P$ and $P$ are
oriented so that $\Pi$ and the inclusion map $i: \tilde P \to \HH$
are orientation preserving.  Cutting $P$ open along the shortest
paths $M_i$'s joining $L_j$ to $L_k$ ($i \neq j \neq k \neq i)$,
we obtain two right-angled hexagons in $P$. Let $Q$ be a lift of
one of the hexagons in $P$ to $\tilde P$ so that  $Q$ is bounded
by complete geodesics $\tilde L_i$ and  $M_i^*$ with $\Pi(\tilde
L_i) =L_i$ and $\Pi(M_i^* \cap Q) = M_i$.
 We choose the lift $Q$ (of one of the
hexagons) so that the cyclic order $\tilde L_1 \rightarrow \tilde
L_2 \rightarrow \tilde L_3$ coincides with the orientation of $Q$.
Let $R_i$ be the hyperbolic reflection about the geodesic $M_i^*$.
Then $\gamma_i = R_{i+2} R_{i+1}$ is the deck transformation group
element so that $\gamma_i(\tilde L_i) = \tilde L_i$ and $\gamma_i$
corresponds to the oriented loop $L_i$.  The closure of the region
in $\HH$ bounded by  $\tilde L_1$, $\tilde L_2$ and $\tilde L_3$
intersects the circle at infinity $S^1_{\infty}$ of $\HH$ in three
disjoint intervals $I_1, I_2, I_3$ where $I_i$ is disjoint from
the closure of $\tilde L_i$. See figure \ref{figure 7}(a). It is known that for
$n \neq 0$
\begin{equation}\label{eqn:gammapower1}
\gamma_i^n(I_i) \subset I_{i+1} \cup I_{i+2}
\end{equation}
and \begin{equation}\label{eqn:gammapower2} \text{end points of
$\gamma_{i+1}^m(\tilde L_i)$ are in $I_i$ for $m>0$}
\end{equation}
where indices are counted modulo 3.  See for instance
\cite{gilman} for a proof.

%We indicate briefly the argument for (\ref{eqn:gammapower1}).
%Assume for simplicity that $i=3$. Let $\tilde M_j$ be the lift of
%$M_j$ so that $\tilde M_j$ appears in the boundary of the hexagon
%$Q$ and let $M_j^*$ be the complete geodesic in $\HH$ which
%contains $\tilde M_j$. Then $I_3$ is in the region $X$ in $\HH
%\cup S^1_{\infty}$ bounded by $M_1^*$ and $M_2^*$. Furthermore, by
%definition of $X$ and $Q$, $X$ is inside a fundamental domain for
%the action of the cyclic group $< \gamma_3
%>$ on $\HH \cup S^1_{\infty} - \{ \text{ end points of $\tilde
%L_3$}\}$. Thus (\ref{eqn:gammapower1}) follows. The statement
%(\ref{eqn:gammapower2}) follows by a similar argument.

\medskip

In the rest of the subsection,  we will focus on $W(L_2, M_3)$
(i.e., i=2, j=3). The general case of $W(L_i, M_j)$ is exactly the
same.

For simplicity, we let $l:=l_2$ and $m:=m_3$. After conjugation by
an isometry of $\HH$, we may assume that $\tilde L_2=G[\infty,
0]$, $\tilde L_3 =G[e,f]$,  $\tilde L_1 =G[c,d]$ with $0<e<f<c<d$, and $\gamma_2(c)=1$.
Since $\gamma_2(z) = e^{-l}z$, we have  $c = e^l$. Note that
$I_1=[0,e]$ in this case. By (\ref{eqn:gammapower2}),
$\gamma_2(\tilde L_1)=G[1, e^{-l}d]$ has end points in $I_1$,
i.e., $1<e^{-l} d < e$. Also, since the distance between $\tilde L_2$
and $\tilde L_1$ is $m$ and $c=e^l$, we  have $d=e^l\coth
^2(\frac{m}{2})$. See figure \ref{figure 7}(b).

\medskip

Define the subset $\Omega_{2,3}$ of $S(\HH)$ as follows. Given $x,
y$ with $0<x<1$ and $c<y<d$, let $q$ be the intersection point
$G[y,x] \cap \tilde L_1$ and let $p$ be the point on $G[y,x]$ so
that the Euclidean ray $0p$ is tangent to the semi-circle
$G[y,x]$. If $\gamma$ is a geodesic path and $v = \gamma'(t)$, we
denote it by $v \in \gamma$.
 Then
\begin{equation}\label{eqn:123} \Omega_{2,3} = \{ v \in S(\HH) | \\ v \in G[y,x], \\
0<x<1,\\ c<y<d, \\ \text{and } \\  \\ v \in
G[q,p]\},\end{equation} where $c=e^l,d=e^l \coth^2(m/2)$.

The main result in this subsection is the following:

\begin{prop}\label{prop:thesetwlm}
  Let $\Pi_* =D \Pi$ be the derivative of the
universal covering map $\Pi: \tilde P \to P$. Then $\Pi_*$ induces
a bijection from $\Omega_{2,3}$ to $W(L_2, M_3)$.  In particular,
the volume of $W(L_2, M_3)$ is $\mu(\Omega_{2,3})$.
\end{prop}

\begin{proof}
We will first show that $\Pi_*(\Omega_{2,3}) \subset W(L_2, L_3)$
and then show that $\Pi_*|_{\Omega_{2,3}}$ is a bijection.

To see $\Pi_*(\Omega_{2,3}) \subset W(L_2, L_3)$, take a vector $v
\in \Omega_{2,3}$ so that $v \in G[q,p] \subset G[y,x]$ as
 in (\ref{eqn:123}).

\begin{lem}\label{lem:betaisasimpleloop}
Let $\tilde \beta =G[ q, \gamma_2( q)]$ be the geodesic in $\HH$
from $ q \in \tilde L_1$ to $\gamma_2( q)$. Then the projection
$\Pi(\tilde \beta) =\beta$ is a simple geodesic loop in $P$ based
at $q' = \Pi(q)$.
\end{lem}

\begin{proof}  By proposition \ref{prop:basicfactssimplegeodesics}, it suffices to show that
 that $\beta \simeq \delta$ rel($q'$) where
$\delta$ is a simple loop at $q'$.
 Indeed, consider the shortest
path $a_1=G[ q,  q_1]$ from $ q$ to $q_1 \in \tilde L_2$. Since
$dist(\tilde L_1, \tilde L_2) = dist(L_1, L_2)$, the projection
$\Pi(a_1)$ is homotopic to $M_3$, the shortest path from $L_2$ to
$L_1$. Thus, by proposition \ref{prop:basicfactssimplegeodesics},
$\Pi(a_1)$ is a simple arc from $L_1$ to $L_2$. Now by the
construction, $\tilde \beta$ and the path $a_1*G[q_1, \gamma_2
(q_1)]* \gamma_2(a_1^{-1})$ have the same end points in $\HH$. Thus
$\beta \simeq \Pi(a_1)*\Pi(G[q_1, \gamma_2 (q_1)])*
\Pi(\gamma_2(a_1)^{-1})$ rel($q'$). Since $\Pi(a_1)$ is an embedded
arc whose interior is disjoint from $\Pi(G[q_1, \gamma_2(q_1)])$
(=$L_2$), by a small perturbation, the loop $\Pi(a_1)*\Pi(G[q_1,
\gamma_2 (q_1)])* \Pi(\gamma_2(a_1)^{-1})$ is relatively homotopic
to a simple loop $\delta$ based at $q'$. It follows that $ \beta
\simeq \delta$ rel($q'$) where $\delta$ is simple. See figure
\ref{figure homotopicloops}. \end{proof}

\begin{figure}[ht!]
\centering
\includegraphics[scale=0.7]{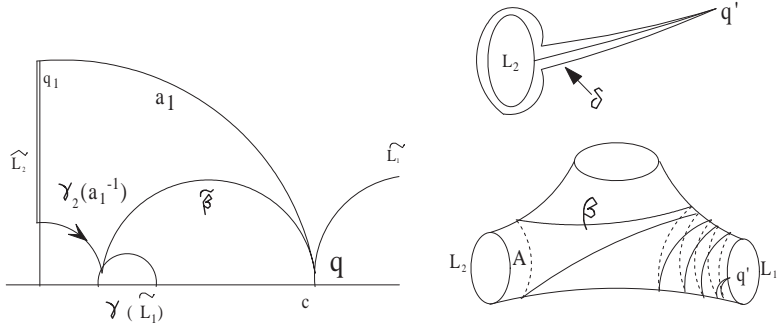}
\caption{homotopic loops are simple} \label{figure homotopicloops}
\end{figure}

Since the simple loop $\beta$ is disjoint from $L_2$ and is
homotopic to $L_2$, there is an annulus $A$ in $P$ bounded by
$\beta$ and $L_2$. Note that $A$ has convex boundary.  The
universal cover $\tilde A$ of $A$ can be identified with the
convex region in $\HH$ bounded by $\tilde L_2$ and the simple path
$\cup_{ n \in \bf Z} \gamma^n_2(\tilde \beta)$.  Now we show that
$\Pi_*(v) \in W(L_2, M_3)$. Consider the geodesic $\gamma(t)
=\Pi_*(G[q, x])$ in $P$ where $\gamma(0)=\Pi(q)$. For $t$ small,
by the construction, $\gamma(t)$ is in the annulus $A$. Since the
vector $v$ is assumed to be generic, $\gamma(t) \in
\partial P$ for some $t>0$. Thus there is the largest $T \in (0, \infty)$  so that
$\gamma([0, T]) \subset A$. First $\gamma(T)$ cannot be in $ L_2$.
Indeed, if this occurs, since $A$ is an annulus, $\gamma|_{[0, T]}
\cong M_3$.  This implies that $G[q,x]$ intersects $\tilde L_2$ and
contradicts $x>0$. It follows that $\gamma(T) \in \beta$. We claim
that, $\gamma_{[0, T]}$ cannot be a simple arc. Otherwise, since $A$
is an annulus and $\gamma|_{[0,T]}$ is an arc joining the same
boundary component $\beta$ of $A$, $\gamma_{[0, T]} \cong \delta$
where $\delta$ is a simple geodesic arc in $\beta$. This contradicts
Gauss-Bonnet theorem since there will be a bi-gon bounded by
$\delta$ and $\gamma|_{[0, T]}$ in the annulus $A$.  It follows that
$\gamma|_{[0, T]}$ is not simple. Let $0<T_1 \leq T$ be the time so
that $\gamma|_{[0, T_1]}$ is a lasso based at $\Pi(q)$ inside $A$.
The loop of this lasso is homotopic to $L_2$ which is the only
simple closed geodesic in $A$. Furthermore, the mid-point of the
loop of the lasso $\gamma|_{[0, T_1]}$ lifts to a point in $G[y,x]$
which is closest to the geodesic $\tilde L_2$. Thus the midpoint of
the lasso is $\Pi(p)$. Finally, the geodesic path $G_P^+(
\Pi_*(v))\cup G_P^-(\Pi_*(v))$ is not homotopic to $M_i^{\pm 1}$ and
to $B^{\pm 1}_1$. Indeed, if otherwise, then a lift of this path
with initial point $q$ will end either on $\tilde L_2$ (homotopic to
$M_3$) or $\tilde L_3$ (homotopic to $M_2$), or $\gamma_2(\tilde
L_1)$, or $\gamma_3(\tilde L_1)$ (homotopic to $B_1$). All these
cases contradict the assumption that
 $0<x<1$.  This shows that $\Pi_*(v) \in W(L_2, M_3)$.

Next, we show that $\Pi_*|$ is  onto.  To see this, take a vector $v
\in W(L_2, M_3)$ so that its graph $G_P(v)^o$ is a lasso based at a
point $q'$ in $L_1$. We claim there is a simple geodesic loop
$\beta$ in $P$ based at $q'$ so that $\beta$ intersects the lasso
$G_P(v)$ only at $q'$ and $\beta$ is freely homotopic to $L_2$.
Indeed,  cutting the surface $P$ open along the lasso $G_P(v)$, we
obtain two convex annuli. See figure \ref{figure 9}(b), (c). One of
the annulus, say $A_1$, contains $L_3$ as a boundary component. Let
$q_1$ and $q_2$ be the preimages of $q'$ in $A_1$ and let $c$ be the
arc in the boundary of $A_1$ joining $q_1$ to $q_2$ so that $c$ is
disjoint from the preimage of $L_1$.

\begin{figure}[ht!]
\centering
\includegraphics[scale=0.65]{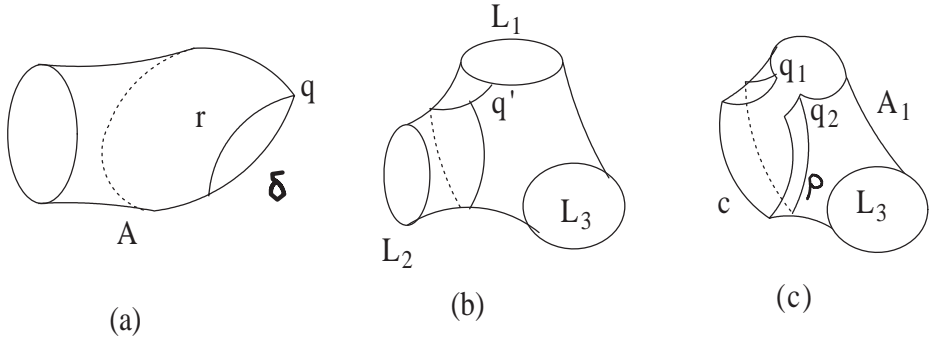}
\caption{cutting surface open along lasso} \label{figure 9}
\end{figure}

 Since $A_1$ is convex, there exists a
shortest geodesic path $\rho$ in $A_1$ joining $q_1$ to $q_2$ so
that $\rho \cong c$ rel($\{q_1, q_2\}$). Since the knot point of
the lasso is a non-smooth point of $c$, the path $\rho$ is
different from $c$. The simple loop $\beta$ is the quotient of
$\rho$ in $P$. Consider $\beta$ as a loop $\beta: S^1 \to P$ and
let $\alpha: \mathbf{R} \to P$ be $\alpha(t) = \beta(e^{it})$.
Since $\beta$ is freely homotopic to $L_2$, there exists a lift
$\tilde \alpha$ of $\alpha$ so that the end points of $\tilde
\alpha$ are the same as that of $\tilde L_2$ in $S^1_{\infty}$.
This lift $\tilde \alpha$ intersects $\tilde L_1$ at exactly one
point $q$ since $\beta \cap L_1 =\{q'\}$. Let $\tilde \gamma(t)$
be the geodesic starting from $q$ which is a lift of the lasso
$G_P(v)$ and $u$ be the unit tangent vector in $\tilde \gamma(t)$
which projects to $v$, i.e., $\Pi_*(u)=v$. We claim that $u \in
\Omega_{2,3}$. Indeed, if $\tilde \beta$ is a lift of the geodesic
path $\beta$ starting at $q$, then $\tilde \alpha$ is the union
$\cup_{n \in \mathbf{Z}} \gamma_2^n (\tilde \beta)$. Let $A$ be
the annulus in $P$ bounded by $\beta$ and $L_2$. Then a universal
cover $\tilde A$ of $A$ is the region bounded by $\tilde L_2$ and
$\cup_{n \in \mathbf{Z}} \gamma_2^n (\tilde \beta)$. It follows
that $\tilde \gamma(t)$ is in $\tilde A$ for $t>0$ small by the
disjointness of $\beta$ and the lasso. Consider the complete
geodesic $G[y,x]$ which contains $\tilde \gamma$. First of all, $
c<y<d$ since $G[y,x]$ intersects $\tilde L_1$. Next, since
$\Pi(G[y,x])$ contains the lasso $G_P(v)$, the preimages of the
knot of $G_P(v)$ in $G[y,x]$ contain two points of the form $z,
\gamma_2(z)$. follows that $\gamma_2(G[y,x]) \cap G[y,x] \neq
\emptyset$, i.e., $ x < e^{-l}y$. But $e^{-l}y < e^{-l}d $. Thus
$x<e^{-l}d$. Next, $x>0$ since $v$ is not in $H(M_i^{\pm 1})$.
Furthermore, it is impossible for $x \in [1, e^{-l}d]$ where $G[1,
e^{-1}d] = \gamma_2(\tilde L_1)$ since $v \notin H(B_i^{\pm 1})$.
Therefore, $0<x<1$. By proposition 4.6, $v$ is between $q'$ and
the midpoint of the lasso. Thus we conclude that $u$ is between
$q$ and $p$. Thus $u \in \Omega_{2,3}$.

Finally, to see that $\Pi_*|$ is injective in $\Omega_{2,3}$,
suppose  that $v_1,  v_2$ in $\Omega_{2,3}$ so that $\Pi_*(v_1)
=\Pi_*(v_2)$. Let $v_i$ be in the geodesic $G[q_i, x_i]$ in
$\Omega_{2,3}$ where $q_i \in \tilde L_1$ and $0<x_i < 1$. Since
$\Pi_*: S(\tilde P) \to S(P)$ is a regular cover with deck
transformation group $\pi_1(P)$, there exists a deck
transformation element $\gamma$ so that $\gamma (v_1) = v_2$. In
particular, $\gamma (G[q_1, x_1]) = G[q_2, x_2]$. This implies
that $\gamma(q_1) = q_2$. Therefore, $\gamma(\tilde L_1) = \tilde
L_1$. However, the only deck transformations leaving $\tilde L_1$
invariant are $\gamma_1^n$. Therefore
$\gamma_1^n(G[q_1,x_1])=G[q_2,x_2]$. If $n \neq 0$, by
(\ref{eqn:gammapower1}) $\gamma_1^n(I_1) \cap I_1 =\emptyset$, we
see that for $x_1 \in (0,1) \subset I_1=[0, e]$, then $x_2 =
\gamma_1^n(x_1) \notin I_1$. Therefore, for $n \neq 0$,
$\gamma_1(G[q_1, x_1])$ cannot be $G[q_2, x_2]$ where $x_1, x_2
\in (0,1)$. This shows that $n=0$, i.e., $v_1=v_2$.
\end{proof}

\subsection{Vectors in  $V(T)$ for a 1-holed torus $T$}\label{ss:vinT}

Let $T$ be a hyperbolic 1-holed torus with geodesic boundary $C$
 and $\{A\}$ the set of non-boundary parallel, simple closed
 geodesics on $T$. Then  $v \in V(T)$ if and only
  if $G_T(v) \cap C \neq \emptyset$. In this case, cutting $T$
   along $G_T(v)^o$ gives a convex hyperbolic cylinder with two
    non-smooth, piecewise geodesic boundaries and there is a unique simple
    closed geodesic $A \subset T$ which is disjoint from $G(v)$.
    Hence $V(T)$ decomposes into the infinite disjoint union
    $V(T)=\bigsqcup_{\{A\}}V_A(T)$ where
$$V_A(T)=\{v \in V(T) ~|~ G_T(v) \cap A =\emptyset\}.$$
Let $P_A$ be the 3-holed sphere obtained by cutting $T$ along $A$
and label the boundaries of $P_A$ so that $L_1=C$, $L_2=A^+$,
$L_3=A^-$. Note that there is an isometric involution of $P_A$
sending $L_2$ to $L_3$ and fixing $L_1$. Then, similar to the
arguments in the previous two subsections, we conclude that
$V_A(T)$ is the disjoint union
$$H(B_1^{\pm 1}) \bigsqcup \bigsqcup_{ i \neq j \neq 1 \neq i} (W(L_i, M_j) \cup W(L_i^{-1},
 M_j))
\bigsqcup \mathbf{A}(\bigsqcup_{ i \neq j \neq 1 \neq i} (W(L_i,
M_j) \cup W(L_i^{-1}, M_j)) ).$$ It follows, from the symmetry of
$P_A$, that
$$\mu(V_A(T))=\mu (H(B_1^{\pm 1}))+ 8\mu(W(L_2,M_3)).$$
Using the notation from \S \ref{ss:pop} and \ref{ss:OHT} that
$length(B_1)=p_1=p_A$, $l_2=l_3=a$ and $m_2=m_3=m_A$, we obtain
$$g(T):= \mu(W(T)) =\mu(S(T))-\mu(V(T)) = 4\pi^2 -\sum_{A}
\mu(V_A(T)).$$ Therefore,
\begin{equation}\label{eqn:gTT}
g(T)=4\pi^2 -8 \sum_{A}( \L(\frac{1}{\cosh^2(p_A/2)}) +
\mu(W(L_2,M_3))).\end{equation}
 This is the formula (\ref{eqn:definitionofg}).
%\end{document}

%% file: sec6.tex
%\section{ Volume calculation}

\section{Calculating the lasso function $La(l,m)$}\label{s:calculatingLa}

%{\bf Tan: Perhaps we need to add an extra line here to explain why we are doing this computation how it is related to $L_i^j$. I think I will leave this section to you so will not make any changes other than the couple of comments below about the Haar measure.}

 By \S \ref{ss:understandingWLM} and the work
 of Bridgeman, we see that the computation of the functions
 $f$ and $g$ reduces to the computation of $\mu(W(L_i,M_j))$
 for a 3-holed sphere $P$, or equivalently, $\mu(\Omega_{i,j})$.
We will show that the volume $\mu(W(L_i, M_j))$ depends only on
the lengths $l_i$, $m_j$ of $L_i$ and $M_j$. The \it lasso
function \rm $La(l_i, m_j)$ is defined to be $\frac{1}{2} \mu(
W(L_i, M_j))$. The goal of this section is to derive an explicit
formula for $La(l,m)$.

  Let us begin by recalling some
  well-known facts about hyperbolic
geometry and the invariant measure on $S(\HH)$.
%
%
%We will use the following notations. We use $\bf H$ to denote the
%upper-half plane model of the hyperbolic geometry. If $x \neq y
%\in $ $\HH \cup
%\partial \HH$, we use $G[x,y]$ to denote the geodesic
%determined from $x$ to $y$. In particular, if $x, y \in \HH$, then
%$G[x,y]$ is the geodesic path from $x$ to $y$.
% The unit tangent
%bundle of $\HH$ will be denoted by $S(\HH)$.
 The invariant measure on the unit tangent bundle
$S(\HH)$ in local coordinates can be written as
$$  \frac{ 2 dx dy du} { (x-y)^2},$$
where $x \neq y \in \bf R$ and $u \in \bf R$  so that the oriented
geodesic $\gamma(v)$ determined by $v \in S(\HH$) is $G[x,y]$ and
$u$ is the signed distance from the base point of $v$ to the
highest point in the semi-circle $G[x,y]$ (in the Euclidean
plane).  See figure \ref{fig:coordinatesofSHH} below.

\begin{figure}[ht!]
\centering
\includegraphics[scale=0.6]{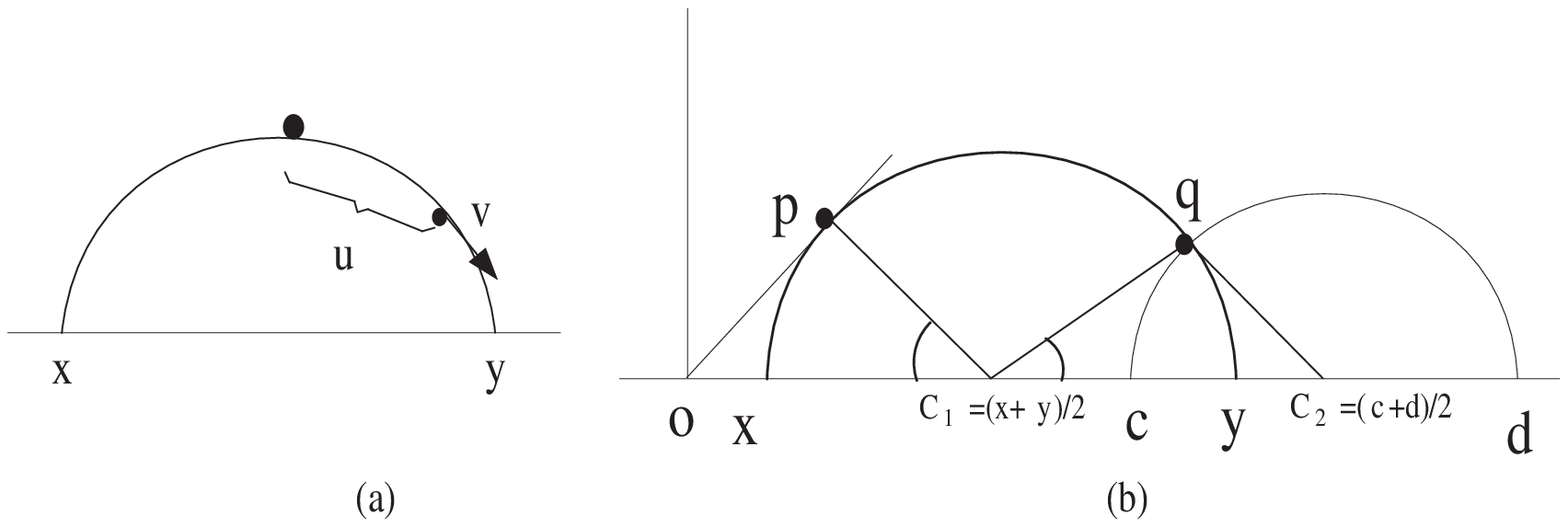}
\caption{Coordinates for $S(\HH)$}\label{fig:coordinatesofSHH}
\end{figure}

Let $\Omega$ be the set defined by (\ref{eqn:123}) (i.e., $\Omega
=\Omega_{2,3}$) where $1<c<d$.
%Given two disjoint geodesics $G[0, \infty]$ and $G[c,d]$ where
%$c,d \in \bf R$ and $1 < c < d$, consider the set of all geodesics
%$G[x,y]$ where $ 0<x<1, c < y< d$. Let $q$ be the point  $\{q\} =
%G[x,y] \cap G[c,d]$ and $p \in G[x,y]$ be the point so that the
%Euclidean ray $Op$ is tangent to $G[x,y]$ at $p$  where $O$ is the
%origin of the Euclidean plane.  Then the set $\Omega$ whose volume
%we are interested in computing is,
%$$ \Omega =\{ v \in S(\HH)  | \text{$\gamma(v) =G[x,y], x \in (0,1), y \in
%(c,d),$}
% \text{the base point of $v$ is in $G[p,q]$} \}.$$
The main result in this section shows
 that the volume $\mu(\Omega)$ of $\Omega$ is
$$\int_0^1 \int_c^d \frac{ \ln | \frac{
y(x-c)(x-d)}{x(y-c)(y-d)}|}{(x-y)^2} dy dx =2 ( \L(\frac{d-1}{d})
- \L(\frac{c-1}{c}) + 2\L(\frac{c-1}{d-1}) -2\L(\frac{c}{c-d}) )$$
where $\L(x)$ is the Roger's dilogarithm.  The right-hand-side of
the above identity will be shown in lemma \ref{lem:simplificationforlassofunction} to be $2[\L(x)
-\L(\frac{1-x}{1-xy}) + \L(\frac{1-y}{1-xy})]$ where
$c=\frac{1}{x}$ and $d=\frac{1}{xy}$. By proposition
\ref{prop:thesetwlm}, we obtain $La(l_i, m_j) =
\frac{1}{2}\mu(W(L_i,M_j)) = \frac{1}{2} \mu(\Omega) =\L(y) -
\L(\frac{1-x}{1-xy}) + \L(\frac{1-y}{1-xy})$ where $c=e^{l_i}$,
$d=e^{l_i}\coth^2(\frac{m_j}{2})$ and $x =e^{-l_i}$,
$y=\tanh^2(m_j/2)$. Combining with lemma \ref{lem:VPdecomposition} and (\ref{eqn:gTT}), we obtain
the formulas  (\ref{eqn:definitionofffirst}) and (\ref{eqn:definitionofg}) in \S\ref{s:functions}.

\subsection{Deriving the volume formula for $\Omega$}

We will establish,

\begin{prop}\label{thm:prop1}
The volume of $\Omega$ is given by

\begin{equation}\label{eqn:VolofOmega}
 \int_0^1 (\int_c^d  \frac{ \ln | \frac{
y(x-c)(x-d)}{x(y-c)(y-d)}|}{(y-x)^2} dy) dx.
\end{equation}
\end{prop}

\begin{proof}   We will use the following known distance formula in the
hyperbolic plane. Namely, $d(e^{i \phi}, e^{-i \psi}) = \ln
\cot(\phi/2) + \ln \cot(\psi/2)$ in $\HH$ where $\phi, \psi \in
(0, \pi/2)$.  Let  $C_1=\frac{x+y}{2}$ be the Euclidean center of
the semi-circle $G[x,y]$ and $\psi$ and $\phi$ be the angles
$\angle 0C_1p$ and $\angle qC_1y$ as shown in figure \ref{fig:coordinatesofSHH}. Then by
the definition of the volume form, we see that $\mu(\Omega)$ is
given by
\begin{equation}\label{eqn:VolofOmega2}
 \int_0^1 (\int_c^d  \frac{ 2\ln \cot(\psi/2) + 2 \ln
\cot(\phi/2)}{(y-x)^2} dy) dx = \int_0^1 (\int_c^d  \frac{ \ln
[\cot^2(\psi/2)\cot^2(\phi/2)]}{(y-x)^2} dy) dx.
\end{equation}

We  calculate $\cot^2(\psi/2)$ and $\cot^2(\phi/2)$ using the
cosine law for Euclidean triangles $\Delta 0pC_1$ and $\Delta
C_1C_2q$ where $C_2=\frac{c+d}{2}$ is the center of the
semi-circle $G[c,d]$.

\medskip
\begin{lem}\label{lem:eulideantrigonometry}
Suppose the lengths of a Euclidean triangle are $l,m, n$ so that
the angle facing the edge of length $l$ is $\theta$. Then
$$ \cot^2(\theta/2) = \frac{ (m+n+l)(m+n-l)}{(m+l-n)(n+l-m)}.$$
\end{lem}

\medskip

Indeed, by the cosine law that $\cos(\theta) = \frac{
m^2+n^2-l^2}{2 mn}$, we obtain
$$ \cot^2(\theta/2) = \frac{ 1+\cos(\theta)}{1-\cos(\theta)}
= \frac{ (m+n)^2-l^2}{ l^2 - (m-n)^2}=\frac{
(m+n+l)(m+n-l)}{(m+l-n)(n+l-m)}.$$

\medskip
For the angle $\psi$, the triangle $\Delta 0pC_1$ is right-angled.
By taking $\theta=\psi$, $n=\frac{y-x}{2}, m =\frac{x+y}{2}$ and
$l=\sqrt{m^2-n^2}$, we obtain $ \cos(\psi) =\frac{y-x}{y+x}$ and
\begin{equation}\label{eqn:cotsquarepsi}
\cot^2(\frac{\psi}{2}) =\frac{y}{x}. \end{equation}

For the angle $\phi$, the edge lengths of the Euclidean triangle
$\Delta q C_1 C_2$ are $n=\frac{y-x}{2}$, $m=\frac{c+d-x-y}{2}$,
and $l=\frac{d-c}{2}$ so that $\phi$ is facing the edge of length
$l$. Now using
$$ l+m+n = \frac{ d-c+y-x+c+d-x-y}{2} = d-x,$$
$$ m+n-l=\frac{ c+d-x-y+ y-x - d +c}{2}= c-x,$$
$$ l+n-m=\frac{d-c+y-x-c-d+x+y}{2}=y-c,$$
$$ l+m-n=\frac{d-c+c+d-x-y -y+x}{2} = d-y,$$
we obtain that
\begin{equation}\label{eqn:cotsquarephi}
\cot^2(\phi/2) =| \frac{ (x-c)(x-d)}{(y-c)(y-d)}|.
\end{equation}
Putting  (\ref{eqn:cotsquarepsi}), (\ref{eqn:cotsquarephi}) into (\ref{eqn:VolofOmega2}), we obtain
\begin{equation}\label{eqn:twologcotpsiovertwo}
 2 \ln \cot(\psi/2) +
 2\ln \cot(\phi/2) = \ln [\cot^2(\psi/2)\cot^2(\phi/2)]
 = \ln | \frac{ y(x-c)(x-d)}{x(y-c)(y-d)}|.
 \end{equation}
 \end{proof}

\subsection{ Evaluation of the integral (\ref{eqn:VolofOmega})}\label{ss:evalofvolofomega}

The evaluation of the integral is similar to the work in
\cite{bridge}. Recall the Roger's dilogarithm $\L$ is defined by
 $\L(0)=0$ and
%\begin{equation}\label{eqn:twodilogderivativex}
$2 \L '(x) = \frac{\ln |x|}{x-1} - \frac{ \ln |x-1|}{x}, \quad
\text{for $x<1$}.$
%\end{equation}
\medskip
\begin{prop}\label{prop:formulafordoubleintegral}
 If $d >c>1$, then
\begin{equation}\label{eqn:formulafordoubleintegralincandd}
 \int_0^1 (\int_c^d \frac{ \ln
|\frac{y(x-c)(x-d)}{x(y-c)(y-d)}|}{(y-x)^2} dy) dx = 2[
\L(\frac{d-1}{d}) - \L(\frac{c-1}{c}) + 2\L(\frac{c-1}{d-1})
-2\L(\frac{c}{c-d})].
\end{equation}
\end{prop}

\begin{proof} To simplify notation, we use
$$R =|\frac{ y(x-c)(x-d)}{x(y-c)(y-d)}|$$ and the integral (\ref{eqn:formulafordoubleintegralincandd}) can
be written as $ \int_0^1\int_c^d \frac{ \ln R}{(x-y)^2} dy dx$.
For simplicity, we drop the constant term in the indefinite
integrals in the lemma below.
\medskip
\begin{lem}\label{lem:integraloflogRoverxminusysquare}
$$\int \frac{\ln R}{(x-y)^2} dy
=\frac{ \ln |\frac{(x-c)(x-d)}{x}|}{x-y} + (\frac{1}{x-y}
-\frac{1}{x}) \ln |y| $$$$+(\frac{1}{x} - \frac{1}{x-c}
-\frac{1}{x-d}) \ln |y-x| +(-\frac{1}{x-y} + \frac{1}{x-c}) \ln
|y-c| +(- \frac{1}{x-y} + \frac{1}{x-d}) \ln |y-d|.$$
\end{lem}

\begin{proof} Using integration by parts, we obtain
$$\int \frac{ \ln R}{(x-y)^2} dy = \int \ln R d(\frac{1}{x-y})$$
$$ = \frac{ \ln R}{x-y} - \int \frac{ d \ln R}{x-y}$$
$$=\frac{ \ln R}{x-y} + \int \frac{dy}{y-x} (\frac{1}{y}
-\frac{1}{y-c} -\frac{1}{y-d})$$ \begin{equation}\label{eqn:logRoverxminusy}
 =\frac{ \ln
R}{x-y} + \int( \frac{1}{(y-x)y} -\frac{1}{(y-x)(y-c)}
-\frac{1}{(y-x)(y-d)}) dy. \end{equation}

Now using the integral formula that for $a \neq b$,
$$ \int \frac{dy}{(y-a)(y-b)} = \frac{1}{(a-b)}[ \ln |y-a| - \ln
|y-b|], $$ we can write (\ref{eqn:logRoverxminusy}) as
$$ \frac{\ln |\frac{(x-c)(x-d)}{x}|}{x-y} +
\frac{ \ln | \frac{y}{(y-c)(y-d)}|}{x-y}
$$
$$ + \frac{1}{x}(\ln |y-x| - \ln |y|) -\frac{1}{x-c}( \ln |y-x| - \ln |y-c|)
-\frac{1}{x-d}( \ln |y-x| - \ln |y-d|)$$
$$=\frac{ \ln |\frac{(x-c)(x-d)}{x}|}{x-y} + (\frac{1}{x-y}
-\frac{1}{x}) \ln |y| $$$$+(\frac{1}{x} - \frac{1}{x-c}
-\frac{1}{x-d}) \ln |y-x| +(-\frac{1}{x-y} + \frac{1}{x-c}) \ln
|y-c| +( -\frac{1}{x-y} + \frac{1}{x-d}) \ln |y-d|.$$ \end{proof}

\medskip
\begin{lem}\label{lem:Wxformula}
 Let $W(x) = \int_c^d \frac{ \ln R}{(x-y)^2} dy$. Then
$$ W(x) = (\frac{ \ln| \frac{x-d}{d}|}{x} - \frac{\ln
|\frac{x}{d}|}{x-d})-  (\frac{ \ln| \frac{x-c}{c}|}{x} - \frac{
\ln |\frac{x}{c}|}{x-c})+  2(\frac{ \ln| \frac{x-c}{x-d}|}{x-d} -
\frac{\ln |\frac{x-d}{x-c}|}{x-c}).$$
\end{lem}

\begin{proof} By lemma \ref{lem:integraloflogRoverxminusysquare}, we can
write $W(x)$ as
$$ (\frac{1}{x-d} - \frac{1}{x-c}) \ln ( | \frac{(x-c)(x-d)}{x}|)
+ (\frac{1}{x-d} - \frac{1}{x}) \ln |d| - ( \frac{1}{x-c} -
\frac{1}{x}) \ln |c|$$
 $$
 +
 (\frac{1}{x} -\frac{1}{x-c} -\frac{1}{x-d})( \ln |x-d| - \ln |x-c|)
+ (-\frac{1}{x-d}+\frac{1}{x-c})\ln |d-c|$$ $$ - \lim_{ y \to c} (
-\frac{1}{x-y} + \frac{1}{x-c}) \ln |y-c|$$
\begin{equation}\label{eqn:limytodminusoneoverxminusy}
 + \lim_{y \to d} (-\frac{1}{x-y} + \frac{1}{x-d}) \ln |y-d|
-(-\frac{1}{x-c} + \frac{1}{x-d}) \ln |c-d|. \end{equation}

Now both limits appearing in (\ref{eqn:limytodminusoneoverxminusy}) are zero since $\lim_{t \to 0}
t \ln |t| =0$. Thus, by rewriting (\ref{eqn:limytodminusoneoverxminusy}) after regrouping according
to $\frac{1}{x}, \frac{1}{x-c}$ and $\frac{1}{x-d}$, we obtain,

$$W(x)=
\frac{1}{x}( -\ln |d| +\ln |c| +\ln |x-d| - \ln |x-c|)
$$
$$ + \frac{1}{x-c}( -\ln |x-c|-\ln |x-d| +\ln |x| -\ln |c| -\ln
|x-d| +\ln |x-c| +\ln |d-c| +\ln |d-c|)$$
$$+ \frac{1}{x-d}( \ln |x-c|+\ln |x-d| -\ln |x| +\ln |d| -\ln
|x-d| +\ln |x-c| -\ln |d-c| -\ln |d-c|)$$
$$ =\frac{1}{x}( \ln |\frac{x-d}{d}| - \ln |\frac{x-c}{c}|)
+\frac{1}{x-c}(\ln |\frac{x}{c}| -2 \ln |\frac{x-d}{c-d}|)
 +\frac{1}{x-d}(-\ln |\frac{x}{d}| +2 \ln |\frac{x-c}{d-c}|)$$
$$=
(\frac{ \ln| \frac{x-d}{d}|}{x} - \frac{\ln |\frac{x}{d}|}{x-d})-
(\frac{ \ln| \frac{x-c}{c}|}{x} - \frac{ \ln |\frac{x}{c}|}{x-c})+
2(\frac{ \ln| \frac{x-c}{x-d}|}{x-d} - \frac{\ln
|\frac{x-d}{x-c}|}{x-c}).
$$
\end{proof}

\medskip

Now to finish the proof of proposition \ref{prop:formulafordoubleintegral}, following
\cite{bridge}, we introduce the following function for $a \neq b$
$$ J(x,a,b) =2 \L(\frac{x-b}{a-b})$$  so that
$$ J'(x,a,b) =\frac{dJ(x,a,b)}{dx} =
\frac{ \ln| \frac{x-b}{a-b}|}{x-a} - \frac{\ln
|\frac{x-a}{b-a}|}{x-b}.$$

By lemma \ref{lem:Wxformula}, it follows that $$W(x) = J'(x,0,d)-J'(x,0,c) +
2J'(x,d,c).$$

Therefore, by the construction the double integral,
$$ \int_0^1 \int_c^d \frac{ \ln R}{(x-y)^2} dy dx $$ $$= \int_0^1 W(x)
dx = J(1,0,d)-J(0,0,d) -J(1,0,c) + J(0,0,c) + 2J(1,d,c) -2J(0,d,c)
$$

But $J(0,0,k) = 2\L(1)$, it follows that

$$\int_0^1 \int_c^d \frac{ \ln R}{(x-y)^2} dy dx =2
( \L(\frac{d-1}{d}) - \L(\frac{c-1}{c}) + 2\L(\frac{c-1}{c-d})
-2\L(\frac{c}{c-d})).$$ \end{proof}

\medskip
To express the volume in terms of the lengths $l$ and $m$,  we
take $c=e^l$ and $d = e^l \coth^2(m/2)$. Then $ \frac{c-1}{c} =
1-e^{-l},$ $ \frac{d-1}{d} = 1- e^{-l} \tanh^2(m/2),$ $
\frac{c-1}{c-d} = (e^{-1}-1) \sinh^2(m/2),$ and $ \frac{c}{c-d} =
-\sinh^2(m/2).$ Thus using
(\ref{eqn:formulafordoubleintegralincandd}), we see the volume
$\mu(\Omega)$ in this case is
$$ 2[\L(   1- e^{-l} \tanh^2(m/2)                   ) -\L(1-e^{-l}            )
+2\L(    (e^{-1}-1) \sinh^2(m/2)                     )
-2\L(-\sinh^2(m/2) )].$$

%This is twice of the function $La(l,m)$ defined in \S
%\ref{ss:identity}, so the total measure of the vectors $v$ generating
% lassos with base point on $L_1$ and loop homotopic to
% $L_2^{\pm 1}$  and such that $v \notin \cup_{l=1}^3
%H(M_l^{\pm 1}) \cup H(B^{\pm 1}_l)$ is
% given by $8La(l,m)=8La(l_2,m_3)$. This, together
% with the results from \S4 about the decomposition
% of $V(P)$ and $V(T)$ completes
%the proof of Theorem \ref{thm:main} by establishing the fomulas for $f(P)$ and $g(T)$ in (\ref{eqn:definitionofffirst}) and (\ref{eqn:definitionofg}).

We now establish the identity (\ref{eqn:lasso}) for the lasso function from
(\ref{eqn:formulafordoubleintegralincandd}).

\begin{lem}\label{lem:simplificationforlassofunction}
Suppose $c=\frac{1}{s}$ and $d=\frac{1}{st}$ in proposition \ref{prop:formulafordoubleintegral}
where $ 1<s,t<1$. Then
\begin{equation}\label{eqn:dilogdminusoneoverd}
\L(\frac{d-1}{d}) - \L(\frac{c-1}{c}) + 2\L(\frac{c-1}{c-d})
-2\L(\frac{c}{c-d}) = \L(t) -\L(\frac{1-s}{1-st}) +
\L(\frac{1-t}{1-st}).
\end{equation}
\end{lem}

\begin{proof}  We have $\frac{d-1}{d} = 1-st$, $\frac{c-1}{c} =
1-s$, $\frac{c-1}{c-d} = \frac{-r}{1-r}$ where $r
=\frac{t(1-s)}{t-1}$ and $\frac{c}{c-d} = -\frac{t}{1-t}$.  Note
 the Roger's dilogarithm satisfies $\L(1-u) = \pi^2/6 -\L(u)$ and
$\L(-\frac{u}{1-u}) = -\L(u)$ for $ 0<u<1$. It follows that
$\L(\frac{c-1}{c-d}) =\L(\frac{-r}{1-r}) =-\L(r)
=-\L(\frac{t(1-s)}{1-st})$ and $\L(\frac{c}{c-d})
=\L(-\frac{t}{1-t}) =-\L(t)$.

Thus the left-hand-side of (\ref{eqn:dilogdminusoneoverd}) is
$$ \L(1-st) - \L(1-s) -
2\L(\frac{t(1-s)}{1-st}) + 2\L(t)$$ $$ =\pi^2/6 -\L(st) -\pi^2/6
+\L(s) -2\L(\frac{t(1-s)}{1-st}) + 2\L(t)$$ $$ =\L(s) -\L(st) +
2\L(t) -2\L(\frac{t(1-s)}{1-st})$$

Using a variation of the pentagon relation (\ref{eqn:pentagon})
that
$$ \L(xy)-\L(x)-\L(y)+\L(\frac{x(1-y)}{1-xy}) +
\L(\frac{y(1-x)}{1-xy})=0,$$ we can write the above as
$$=\L(t) +\L(\frac{s(1-t)}{1-st}) -\L(\frac{t(1-s)}{1-st}).$$
Since $\frac{s(1-t)}{1-st}=1-\frac{1-s}{1-st}$ and
$\L(1-u)=\pi^2/6 -\L(u)$, the above equation is $\L(t)
-\L(\frac{1-s}{1-st}) + \L(\frac{1-t}{1-st})$.
\end{proof}

\begin{cor}\label{cor:eqnforfP}(Equation (\ref{eqn:definitionoff}) for $f(P)$)\label{cor:simplificationforfP}
 Suppose $P$ is a hyperbolic 3-holed sphere of boundary
lengths $l_i$'s so that the lengths of $M_i$ are $m_i$ and $B_i$ are
$p_i$. Let $x_i = e^{-l_i}$ and $y_i=\tanh^2(m_i/2)$. Then
$$ f(P) =8 [\sum_{ i \neq j} (\L(\frac{1-x_i}{1-x_i y_j})
-\L(\frac{1-y_j}{1-x_i y_j}))
-\sum_{k=1}^3(\L(y_k)+\L(\frac{1}{\cosh^2(p_k/2)}))]$$ $$ =4
\sum_{i \neq j}[2\L(\frac{1-x_i}{1-x_i y_j})
-2\L(\frac{1-y_j}{1-x_i y_j})-\L(y_j)
-\L(\frac{(1-y_j)^2 x_i}{(1-x_i)^2y_j})]$$
\end{cor}

\begin{proof}
Recall that by definition and lemma \ref{lem:VPdecomposition}, $f(P) = \mu(W(P)) =
\mu(S(P)) - \mu(V(P)) = 4\pi^2 -[\sum_{i=1}^3 (\mu(H(M^{\pm 1}) +
\mu(H(B^{\pm 1})) + 4\sum_{i \neq j} \mu(W(L_i, M_j)]$.  It
follows that
$$f(P) = 4\pi^2 -8[\sum_{i=1}^3 (\L(\frac{1}{\cosh^2(m_i/2)})
+\L(\frac{1}{\cosh^2(p_i/2)})) +\sum_{ i \neq j} La(l_i, m_j)].$$
Using $\L(\frac{1}{\cosh^2(m_i/2)}) =\L(1-y_i) =\pi^2/6 -\L(y_i)$
and lemma \ref{lem:simplificationforlassofunction}, we can write
the above as
$$ 4\pi^2-8[\sum_{i=1}^3 (\L(1-y_i) +\L(\frac{1}{\cosh^2(p_i/2)}  )) +\sum_{i \neq j}
(\L(y_i) -\L(\frac{1-x_i}{1-x_i y_j})
+\L(\frac{1-y_j}{1-x_iy_j}))].$$ Since $\L(1-y_i)+\L(y_i)
=\pi^2/6$, and $\frac{1}{\cosh^2(p_k/2)} =
\frac{(1-y_j)^2 x_i}{(1-x_i)^2y_j}$ by (\ref{eqn:perpendicular}), the above equation is
equivalent to the identity in the corollary.  \end{proof}

 \vskip 30pt

 \subsection{Surfaces with boundary}\label{ss:surfaceswithboundary}
%{This is the revised version 7 Apr 2011}

Let $F$ be a hyperbolic surface with non-empty geodesic boundary such that the Euler characteristic $\chi(F)<-1$.
As in \S 4, for a generic unit tangent vector $v \in S(F)$, $G(v)$ is a graph contained in an embedded simple geometric subsurface $\Sigma$ of $F$, except that now, $G(v) \cap \partial F$ may not be empty. Again, as in \S 4, we need to calculate $\mu(W(\Sigma))$, where
$$W(\Sigma)=\{v \in S(\Sigma)|G(v)=G_{\Sigma}(v)\}.$$
When $\Sigma \cap \partial F =\emptyset$, then the computation of $\mu(W(\Sigma))$ is exactly the same as in \S 4. This occurs when $\Sigma$ is a 1-holed torus $T$ (since $\chi(F)<-1$), or when it is a 3-holed sphere $P$ for which $\partial P \cap \partial F =\emptyset$. It remains to compute $\mu(W(P))$ when $P$ is an embedded geometric 3-holed sphere for which $\partial P \cap \partial F$ consists of either one or two components.

Let $L_1,L_2,L_3$ be the boundary components of $P$. We first consider the case where  $\partial P \cap \partial F$ has one component, which we may take to be $L_1$ by convention. We also use the shorthand notation
\begin{equation}\label{eqn:shorthandwLi}
W(L_i^{ \pm 1},M_j)=W(L_i,M_j)\bigsqcup W(L_i^{ - 1},M_j).
\end{equation}
We see  from the definition of $W(P)$ that in this case, besides  spines $G(v)$ for $P$ which do not intersect $\partial P$, $G(v)$ is also a spine for $P$ when
$$v \in H(B_1^{\pm 1}) \cup  W(L_2^{ \pm 1},M_3) \cup W(L_3^{ \pm 1},M_2) \cup {\mathbf A}(W(L_2^{ \pm 1},M_3) \cup W(L_3^{ \pm 1},M_2)).$$

It follows that for such $P$,
\begin{equation}\label{eqn:hatf2}
    {\hat f}(P) := \mu(W(P))=f(P)+8\big(\L(\frac{1}{\cosh^2 p_1/2})+La(l_2,m_3)+La(l_3,m_2)\big)
\end{equation}

The remaining case is when $\partial P \cap \partial F$ has two components, which we may take to be $L_1$ and  $L_2$ by convention. Now, besides spines $G(v)$ for $P$ which do not intersect $\partial P$, $G(v)$ is also a spine for $P$ when
\begin{eqnarray*}
v &\in& H(B_1^{\pm 1}) \cup H(B_2^{\pm 1}) \cup H(M_3) \\
  && \quad \cup W(L_2^{\pm 1},M_3) \cup W(L_3^{\pm 1},M_2)\cup W(L_1^{\pm 1},M_3) \cup W(L_3^{\pm 1},M_1) \\
 && \quad  \cup {\mathbf A}( W(L_2^{\pm 1},M_3) \cup W(L_3^{\pm 1},M_2)\cup W(L_1^{\pm 1},M_3) \cup W(L_3^{\pm 1},M_1)).
 \end{eqnarray*}

It follows that for such $P$,

\begin{eqnarray*}\label{eqn:hatf}
    {\bar f}(P) &:=& \mu(W(P)) = f(P)+8\{\L(\frac{1}{\cosh^2 p_1/2})+\L(\frac{1}{\cosh^2 p_2/2})+\L(\frac{1}{\cosh^2 m_3/2})\\
    && \qquad \qquad \qquad +La(l_2,m_3)+La(l_3,m_2)+La(l_1,m_3)+La(l_3,m_1)\}
\end{eqnarray*}

Theorem 1.2 now follows.

%% file: ref.tex
{}

%% file: appendix.tex
\bigskip

\bigskip

\noindent  {\bf Appendix}

\medskip

\noindent {\bf A1. Pentagon relations for dilogarithm and
hyperbolic pentagons}\label{ss:pentagonrelations}
The following simple property was discovered during our study of
the Roger's dilogarithm. It puts the pentagon relations in the
perspective of lengths of hyperbolic right-angled pentagons.

\begin{prop}\label{prop:hyperbolicpentagonrelations}
Suppose $l_1, ..., l_5$ are the lengths of the five sides of a
hyperbolic right-angled pentagon. Then

$$ \sum_{i=1}^5 \L(\tanh^2(l_i)) =\pi^2/2,$$
and
$$ \sum_{i=1}^5 \L(\frac{1}{\cosh^2(l_i)}) =\pi^2/3.$$

In fact, each of the above is equivalent to the pentagon relation (\ref{eqn:pentagon})
for the Roger's dilogarithm.
\end{prop}

\begin{proof} We assume all pentagons are right-angled in the sequel.
We begin with the sine law for pentagons. Suppose the edges $e_1,
..., e_5$ in the hyperbolic pentagon are cyclically labelled so
that the length of $e_i$ is $l_i$. Then the sine law for pentagon
says that $\cosh^2(l_i) = \sinh(l_{i+2})\sinh(l_{i+3})$. Let $s_i
=\sinh^2(l_i)$, then the sine law says  $s_{i}+1 =
s_{i+2}s_{i+3}$.  Let $x_i =\tanh^2(l_i)=\frac{s_i}{s_i+1}$ and
let $x=x_1$, $y=x_3$. Then $x_5=\frac{1-x}{1-xy}$,
$x_4=\frac{1-y}{1-xy}$ and $x_2=1-xy$ by the relations
$s_i+1=s_{i+2}s_{i+3}$.

Now the pentagon relation for the Roger's dilogarithm $\L(t)$ says
for $x,y \in (0,1)$,  $$\L(x)+\L(y) +\L(1-xy) +
\L(\frac{1-x}{1-xy}) +\L(\frac{1-y}{1-xy})
 =\frac{\pi^2}{2}.$$  The five variables inside $\L(t)$ are
 exactly $\tanh^2(l_i)$ by the above calculation. Thus the first
 identity follows. Since $\L(\frac{1}{\cosh^2(x)}) =\pi^2/6 -
 \L(\tanh^2(x))$, the second equation follows.

\end{proof}

\bigskip

\noindent {\bf A2. Using different rules to generate $G(v)$}\label{ss:differentrules}

%{This is the revised version 7 Apr 2011}

For a generic unit tangent vector $v \in S(F)$, we gave a somewhat arbitrary rule to define the graph $G(v)$ in \S3 (generating the geodesic at equal speed in both forwards and backwards direction until we obtain intersections), from which we obtained the decomposition of the unit tangent bundle $S(F)$ which gave rise to the identities in Theorems \ref{thm:main} and \ref{thm:mainwithboundary}. The main advantage of our choice was that for generic vectors $v \in S(F)$,  $G(v)=G(-v)$ so that in the computation of the measures $\mu(W(\Sigma))$  for geometrically  embedded simple surfaces in \S\ref{s:measureofdecompostion}, we were able to exploit the symmetry in our computations. In particular, in the computation of the measure of the set of vectors $v \in S(\Sigma)$ which generated  lassos, we just doubled the measure of the vectors which generated the positively oriented  lassos. A natural question which arises is whether we get different identities if we use a different rule for generating $G(v)$. As an example, a fairly natural choice would be a forward first rule, that is, to generate $g^+(v)$ until the first point of intersection, after which we generate $g^{-}(v)$ until the next point of intersection, thereby producing a graph $G(v)$ as in \S\ref{s:decomposing}. More generally, we may generate $g^{+}(v)$ and $g^{-}(v)$ at different fixed constant speeds to obtain $G(v)$.

It is clear that the homotopy type of $G(v)$ may be different for different rules, hence, we would obtain a different decomposition of the unit tangent bundle $S(F)$. We claim here that nonetheless, the resulting identities obtained are all the same. The main observation is that the measure of the complementary set $V(\Sigma)$ of vectors which do not generate spines for a simple surface $\Sigma \subset F$  are the same, for different rules.

We give a brief explanation here.
Recall  from Lemma \ref{lem:VPdecomposition} that $v \in V(P)$ if $v \in H(M_i)$ or $H(B_i)$, $i=1,2,3$, or $v ~~{\hbox{or}}~~-v \in W(L_i^{\pm 1},M_j)$, $1 \le i \neq j \le 3$. Furthermore, the sets are disjoint.

There is no problem with $H(M_i)$ and $H(B_i)$, the sets are the same whatever rules we use to define $G(v)$ and so they have the same measures. The issue arises in the sets $W(L_i^{\pm 1},M_j)$, $1 \le i \neq j \le 3$, which depend on the rule used to define $G(v)$. More specifically, suppose  that $\alpha:[T_1,T_2] \rightarrow P$ is a positively oriented lasso on $P$ with base point on $L_1$ and with a positive loop around $L_2$ such that $\alpha(T_3)=\alpha(T_2)$ for some $T_1<T_3<T_2$ (cf definition \ref{def:lassos}). Then, if we use the original rule for generating $G(v)$, $v=\alpha'(t) \in W(L_2,M_3)$ (that is,  $G(v)=\alpha$) if and only if $T_1<t<\frac{T_3+T_2}{2}$ and $-v$ generates $-\alpha$ if and only if $v$ generates $\alpha$. However, for example, if we use the forward first rule instead, than $v=\alpha'(t) \in W(L_2,M_3)$ (that is,  $G(v)=\alpha$) if and only if $T_1<t<T_3$, while $-v$ generates $-\alpha$ if and only if $v=\alpha'(t)$ with $T_1<t< T_2$.

The main observation is that when we sum over the measure of all $v$ and $-v$ which generate either $\alpha$ or $-\alpha$, it is given by $(T_3+T_2)-2T_1$, which is the same as for the first rule. If we let $W^-(L_2,M_3)$ be the set of vectors $v \in S(P)$ generating  lassos with negative orientation and base point at $L_1$ and loop homotopic to $L_2$, then it follows that the measure of $W(L_2,M_3) \cup W^-(L_2,M_3)$  is the same for both rules. In fact, the same argument shows that any consistently applied rule gives the same measure for $W(L_2,M_3) \cup W^-(L_2,M_2)$, the extra measure in one set is compensated by the deficit in the other. It follows that $f(P)$ depends only on the lengths $L_1,L_2$ and $L_3$. A similar argument holds for $g(T)$, $\bar{f}(P)$ and $\hat{f}(P)$.